\documentclass[11pt, twoside]{article}
\usepackage{latexsym}
\usepackage{amsmath}
\usepackage{amssymb}
\usepackage[all]{xy}
\usepackage{amsfonts}
\usepackage{verbatim}
\usepackage{amsthm}
\usepackage{bm}
\usepackage{mathrsfs}
\usepackage{epsfig}
\usepackage{xy}
\usepackage{array}
\usepackage{stmaryrd}
\usepackage{graphicx,color}
\usepackage{xcolor}
\usepackage{tikz}
\usepackage{xcolor}
\usepackage[colorlinks=true,linkcolor=blue,citecolor=blue]{hyperref}
\usetikzlibrary{arrows,calc}
\usepackage{etex}
\usepackage{mathdots}
\usepackage{float}
\usepackage{graphics}
\usepackage{pdflscape}
\usepackage{CJK}
\usepackage{anysize,hyperref}
\input xypic
\xyoption{all}
\usepackage{extarrows}
\usepackage[perpage,symbol]{footmisc}
\usepackage{setspace}
\usepackage{enumitem}
\setenumerate[1]{itemsep=0pt,partopsep=0pt,parsep=\parskip,topsep=5pt}
\setitemize[1]{itemsep=0pt,partopsep=0pt,parsep=\parskip,topsep=5pt}
\setdescription{itemsep=0pt,partopsep=0pt,parsep=\parskip,topsep=5pt}

\renewcommand{\cal}{\mathcal}

\def\T{\mathcal{T}}

\def\P{\mathcal{P}}
\def\I{\mathcal{I}}

\def\R{\mathcal{R}}
\def\M{\mathcal{M}}
\def\T{\mathcal{T}}
\def\C{\mathscr{C}}

\def\dr{\ar@{->}[r]}

\usepackage{pgf,tikz}

  \usetikzlibrary{arrows}
  \usetikzlibrary{decorations.markings}
  \usetikzlibrary{shapes}
  \usetikzlibrary{matrix}
  \usepackage{wasysym}
\begin{document}
\baselineskip=15pt
\title{\Large{\bf Abelian quotients of categories of $n$-exangles \footnotetext{$^\ast$Corresponding author~~ }}}
\medskip
\author{Yutong $\rm{Zhou^{\ast}}$ }

\date{}

\maketitle
\def\blue{\color{blue}}
\def\red{\color{red}}

\newtheorem{theorem}{Theorem}[section]
\newtheorem{lemma}[theorem]{Lemma}
\newtheorem{corollary}[theorem]{Corollary}
\newtheorem{proposition}[theorem]{Proposition}
\newtheorem{conjecture}{Conjecture}
\theoremstyle{definition}
\newtheorem{definition}[theorem]{Definition}
\newtheorem{question}[theorem]{Question}
\newtheorem{remark}[theorem]{Remark}
\newtheorem{remark*}[]{Remark}
\newtheorem{example}[theorem]{Example}
\newtheorem{example*}[]{Example}
\newtheorem{condition}[theorem]{Condition}
\newtheorem{condition*}[]{Condition}
\newtheorem{construction}[theorem]{Construction}
\newtheorem{construction*}[]{Construction}

\newtheorem{assumption}[theorem]{Assumption}
\newtheorem{assumption*}[]{Assumption}

\baselineskip=17pt
\parindent=0.5cm

\begin{abstract}
\baselineskip=16pt

 The notion of $n$-exangulated categories was introduced by Herschend-Liu-Nakaoka, which is a simultaneous generalization of $n$-exact categories in the sense of Jasso and  $(n+2)$-angulated categories in the sense of Geiss-Kelier-Oppermann. Let $(\C,\mathbb{E},\mathfrak{s})$ be an $n$-exangulated category with enough
projectives $\P$ and $\M$ a full subcategory of $\C$ containing $\P$. In the present paper, It is proved that a certian quotient category of $\mathfrak{s}$-def$(\M)$ is abelian. We denoted by $S(\C)$ the category of $n$-exangles, whose object are given by distinguished $n$-exangles in $\C$. If $\M=\C$, we obtain that a certain ideal quotient category $S(\C)/\R_2$ is equivalent to the category of finitely presented modules mod-$(\C/[\P])$. Furthermore, we present the quotient category $S(\C)/\R_2$ always has an abelian structure when taking $n$ as an even number. The abelian quotient $S(\C)/\R_2$ admits some nice properties. We describe the projective objects in $S(\C)/\R_2$ and characterize the simple objects in $S(\C)/\R_2$ as Auslander-Reiten $n$-exangle sequences in $\C$.\\[0.5cm]
\medskip
\end{abstract}

\pagestyle{myheadings}
\markboth{\rightline {\scriptsize Y. $Zhou^{\ast}$} }
         {\leftline{\scriptsize  Abelian quotients of categories of $n$-exangles }}

\section{Introduction}

 ~~Geiss, Keller and Oppermann \cite{GKO} introduced the notion of an $(n+2)$-angulated category, which is a higher dimensional analogue of a triangulated category and showed that certain $n$-cluster tilting subcategories of triangulated categories are $(n+2)$-angulated. In a more recent development, the concepts of $n$-exact categories and algebraic $n$-angulated categories are introduced by Jasso \cite{j}. He showed that the quotient category of a Frobenius $n$-exact category has a natural structure of an $(n+2)$-angulated category. Moreover, Jasso proved that $n$-cluster-tilting subcategories of abelian(exact) categories are $n$-abelian($n$-exact). Recently, Herschend--Liu--Nakaoka \cite{HLN} defined $n$-exangulated categories as a higher dimensional analogue of extriangulated categories. In particular, we know that $n$-exangulated categories simulataneously generalization of
  $n$-exact categories and $(n+2)$-angulated categories. Hence many results hold on $n$-exact categories and $(n+2)$-angulated categories can be unified in the same framework. In addition, we note that there are some other examples of $n$-exangulated categories which are neither $n$-exact nor $(n+2)$-angulated, for more details,
 we refer to \cite{ HHZ, HLN, HLN1, LZ}.

 Cluster tilting theory provides a way of constructing abelian categories from some triangulated categories, in general,  one can pass from triangulated categories to abelian categories by factoring out cluster tilting subcategories. Let $\C$ be a triangulated category and $\T$ be a cluster-tilting subcategory of $\C$, then the quotient $\C/[\T]$ is abelian, see \cite{BMR, KR, KZ} for more details. Cluster tilting theory is also suitable for exact category. Demonet-Liu \cite{DL} provided a general framework for pass from exact categories to abelian categories by factoring out cluster tilting subcategories. Submodule categories provide another way to construct abelian quotient categories. Certain quotients of submodule categories are realized as categories of finitely presented modules over stabel Auslander algebras \cite{E,RZ}.  Lin \cite{L} showed that some quotients of categories of short exact sequences in exact categories are abelian.  The version of triangulated categories is given in \cite{N}, certain quotients of categories of triangles are abelian.
 More precisely, Lin \cite{Li} investigate abelian quotients arising from extriangulated categories via morphism categories.

 Inspired by above work, in the present paper, let $(\C,\mathbb{E},\mathfrak{s})$ be an $n$-exangulated category, we mainly studies the abelian quotients of categories of $n$-exangles. We denote by $S(\C)$ the category of $n$-exangles, where the objects are the distinguished $n$-exangles in $\C$

$$A^\bullet=(\xymatrix{
 A_0\ar[r]^{\alpha_0} & A_1 \ar[r]^{\alpha_1}& A_2 \ar[r]^{\alpha_2} & \cdots \ar[r]^{\alpha_{n-2}}&A_{n-1} \ar[r]^{\alpha_{n-1}}& A_n \ar[r]^{\alpha_{n}} & A_{n+1} \ar@{-->}[r]^{\delta}&})$$
and the morphisms from $A^\bullet$ to $B^\bullet$ is defined to be the $\varphi\bullet=(\varphi_0,\varphi_1,\cdots,\varphi_{n+1})$ such that the following diagram is commutative
$$\xymatrix{
A_0 \ar[r]^{\alpha_0}\ar[d]^{\varphi_0} & A_1 \ar[r]^{\alpha_1}\ar[d]^{\varphi_1} & A_2 \ar[r]^{\alpha_2}\ar[d]^{\varphi_2} & \cdots \ar[r]^{\alpha_{n-1}}& A_{n} \ar[r]^{\alpha_{n}}\ar[d]^{\varphi_{n}} & A_{n+1} \ar@{-->}[r]^{\delta}\ar[d]^{\varphi_{n+1}} &{~}\\
B_0 \ar[r]^{\beta_0} & B_1 \ar[r]^{\beta_1} & B_2 \ar[r]^{\beta_2} & \cdots \ar[r]^{\beta_{n-1}}& B_{n} \ar[r]^{\beta_{n}} & B_{n+1} \ar@{-->}[r]^{\delta'} & {~}}$$
and ${\varphi_0}_{\ast}\delta={\varphi_{n+1}}^{\ast}\delta'$. Let $A^\bullet$ and $B^\bullet$ be two distinguished $n$-exangles, we denote by $\R_2(A^\bullet,B^\bullet)$ (resp. $\R_1'(A^\bullet,B^\bullet)$) the class of morphism $\varphi^\bullet:A^\bullet\rightarrow B^\bullet$ such that $\varphi_{n+1}$ factors through $\beta_{n}$ (resp. $\varphi_{0}$ factors through $\alpha_{0}$).  Assume that $\M$ is a full subcategory of $\C$. We denote by $\mathfrak{s}$-def$(\M)$ (resp. $\mathfrak{s}$-inf$(\M)$) the full subcategory of Mor($\C$) consisting of $\mathfrak{s}$-deflations (resp.$\mathfrak{s}$-inflations). The full subcategory of $\mathfrak{s}$-def$(\M)$ consisting of split epimorphisms (resp. split monomorphisms) is denoted by s-epi($\M$) (resp. s-mono($\M$)). We denote by sp-epi($\M$)(resp. si-mono$(\M)$) the full subcategory of $\mathfrak{s}$-def$(\M)$ consisting of $(\xymatrix{M \ar[r]^-{1} & M})\oplus(P\rightarrow M')$ (resp. ($\xymatrix{M \ar[r]^-{1} & M})\oplus(M'\oplus I)$) with $P$ projective (resp. $I$ injective). The following is our main theorem.\\
{\bf Theorem 1.1} Let $(\C,\mathbb{E},\mathfrak{s})$ be an $n$-exangulated category.
\begin{itemize}
\item[$\rm{(i)}$] If $\C$ has enough projectives, then we have the following equivalences: $$S(\C)/\R_2\cong\mathfrak{s}\textnormal{-def}(\C)/[\textnormal{s-epi}(\C)]\cong
    \textnormal{mod-}(\C/[\P]).$$

\item[$\rm{(ii)}$] If $\C$ has enough injectives $\I$, then we have the following equivalences: $$S(\C)/\R_2\cong\mathfrak{s}\textnormal{-inf}(\C)/[\textnormal{si-mono}(\C)]
    \cong{(\textnormal{mod-}{{(\C/[\I])}^{op}})}^{op}.$$
\end{itemize}
 {\bf Theorem 1.2} Let $(\C,\mathbb{E},\mathfrak{s})$ be an $n$-exangulated category. When taking $n$ as an even number, the quotient category $S(\C)/\R_2$ is an abelian category.

 This paper is organized as follows.

 In section \ref{2}, we recall some definitions and basic facts of $n$-exangulated categories and morphism categories. Then in section \ref{3}, let $(\C,\mathbb{E},\mathfrak{s})$ be an $n$-exangulated category, we fix some notations in $\C$, the main result of this section is Theorem \ref{th34}, we show that some quotient categories of $\mathfrak{s}$-deflation categories are equivalent to module categories. Finally, in section \ref{4}, we study abelian quotients of the categories of $n$-exangles and then we prove our main result, see Theorem \ref{th44}.

Throughout this paper, We assume that $n$ is an positive integer. Unless otherwise stated, that all subcategories of additive categories considered are full, closed under isomorphisms and direct summands, all functors between additive categories are hypothesised to be additive. The same convention is used for morphisms in an arbitrary category. For a category $\C$, we will also write $\C(X,Y)$ for the set of morphisms between two object $X$ and $Y$ of $\C$.

\section{Definitions and Preliminaries}\label{2}

In this section, we first give some facts on morphism categories.  Then recall some basic results and notions on $n$-exangulated categories in \cite{HLN}.

\subsection{Morphism Categories}
Assume that $\C$ is an additive category, The \emph{morphism category} of $\C$ is the category ${\rm{Mor}}(\C)$ defined by the following statements. The objects of ${\rm{Mor}}(\C)$ are all the morphisms $f:X\rightarrow Y$ in $\C$. The morphisms from $f:X\rightarrow Y$ to $f':X'\rightarrow Y'$ are pairs $(a,b)$ where $a:X\rightarrow X'$ and $b:Y\rightarrow Y'$ such that $bf=f'a$. The composition of morphisms is componentwise.
Let $\C$ be an additive category. A right \emph{$\C$-module} is a contravariantly additive functor $F:\C\rightarrow{\rm{Ab}}$ where $\rm{Ab}$ is the category of abelian groups. Denote by Mod-$\C$ the category of right $\C$-modules. It is well known that Mod-$\C$ is an abelian category. A $\C$-module $F$ is called \emph{finitely presented} if there exists an exact sequence $\C(-,X)\rightarrow\C(-,Y)\rightarrow F\rightarrow0$. We denote by mod-$\C$ the full subcategory of ${\rm{Mor}}(\C)$ formed by finitely presented $\C$-modules, by proj-$\C$ (resp. inj-$\C$) the full subcategory of mod-$\C$ consisting of projective (resp. injective) objects. It is known that  mod-$\C$ is abelian if and only if $\C$ admits weak kernels.

\subsection{$n$-Exangulated Categories}

Let $\C$ be an additive category and $\mathbb{E}:{\C}^{op}\times\C\rightarrow \rm{Ab}$ ($\rm{Ab}$ is the category of abelian groups) an additive bifunctor. For any pair of objects $A,C\in\C$, an element $\delta\in\mathbb{E}(C,A)$ is called an \emph{$\mathbb{E}$-extension} or simply an \emph{extention}. We also write such $\delta$ as $_A{\delta}_{C}$ when we indicate $A$ and $C$.
The zero element $_A{0}_{C}=0\in\mathbb{E}(C,A)$ is called the \emph{{split} $\mathbb{E}$-extension}. Since $\mathbb{E}$ is a bifunctor, for any morphism $a\in\C(A,A^{'})$ and $c\in\C(C^{'},C)$, we denote the $\mathbb{E}$-extension $\mathbb{E}(C,a)(\delta)\in\mathbb{E}(C,A')$ by $a_{\ast}\delta$ and denote the $\mathbb{E}$-extension $\mathbb{E}(c,A)(\delta)\in\mathbb{E}(C',A)$ by $c^{\ast}\delta$. In this terminology, we have
$$\mathbb{E}(c,a)(\delta)=c^\ast a_\ast\delta=a_\ast c^\ast\delta\in\mathbb{E}(C',A').$$
 Let $_A\delta_C$ and $_{A'}\delta'_{C'}$ be any pair of $\mathbb{E}$-extensions. A \emph{morphism}
$(a,c):\delta\rightarrow\delta'$ of $\mathbb{E}$-extensions is a pair of morphism $a\in\C(A,A^{'})$ and $c\in\C(C,C')$ in $\C$, satisfying the equality $$a_{\ast}\delta=c^{\ast}\delta'.$$
Let $_A{\delta}_{C}$ and $_{A'}{\delta}_{C'}$ be any pair of $\mathbb{E}$-extensions, $\delta\oplus\delta'\in\mathbb{E}(C\oplus C', A\oplus A')$ be the element corresponding to $(\delta,0,0,\delta')$ through the natural isomorphism $\mathbb{E}(C\oplus C', A\oplus A')\simeq\mathbb{E}(C,A)\oplus\mathbb{E}(C,A')\oplus\mathbb{E}(C',A)\oplus\mathbb{E}(C',A').$

\begin{definition}\cite[Definition 2.7]{HLN}
Let ${\bf{C}}_{\C}$ be the category of complexes in $\C$. As its full subcategory, define ${\bf{C}}^{n+2}_\C$ to be the category of complexes in $\C$ whose components are zero in the degrees outside of $\{0,1,\cdots,n+1\}$. Namely, an object in ${\bf{C}}^{n+2}_\C$ is a complex $X^\bullet=\{X_i,d^X_i\}$ of the form
$$X_0 \xrightarrow{~d^X_0~} X_1 \xrightarrow{~d^X_1~} \cdots \xrightarrow{~d^X_{n-2}~} X_{n-1} \xrightarrow{~d^X_{n-1}~} X_{n}\xrightarrow{~d^X_n~} X_{n+1}.$$

We write a morphism $f^\bullet:X^\bullet\rightarrow Y^\bullet$ simply as $f^\bullet=(f_0,f_1,\cdots,f_{n+1})$, only indicating the terms of degrees $0,1,\cdots,n+1$.

\end{definition}

\begin{definition}\cite[Definition 2.9]{HLN}
Let $\C,\mathbb{E},n$ be as before. Define a category ${\AE}:={\AE}^{n+2}_{(\C,\mathbb{E})}$ as follows.
\begin{itemize}
\item[(a)] An object in ${\AE}^{n+2}_{(\C,\mathbb{E})}$ is a pair $\langle X^\bullet,\delta\rangle$ of $X^\bullet\in{\bf{C}}^{n+2}_{\C}$ and $\delta\in\mathbb{E}(X_{n+1},X_0)$ is called an  \emph{$\mathbb{E}$-attached complex of lenth $n+2$}, if it satisfies
    $$(d^X_0)_\ast\delta=0~~\text{and}~~(d^X_n)^\ast\delta=0.$$
    We also denote it by
$$\xymatrix{
X_0 \ar[r]^-{d^X_0} & X_1 \ar[r]^-{d^X_1} & \cdots \ar[r]^-{d^X_{n-2}}&X_{n-1} \ar[r]^-{d^X_{n-1}}& X_n \ar[r]^-{d^X_{n}} & X_{n+1} \ar@{-->}[r]^-{\delta}&}$$
\item[(b)] For such pair $\langle X^\bullet,\delta\rangle$ and $\langle Y^\bullet,\rho\rangle$,
\emph{a morphism} $f^\bullet:\langle X^\bullet,\delta\rangle\rightarrow\langle Y^\bullet,\rho\rangle$
is defined to be a morphism $f^\bullet\in{\bf{C}}^{n+2}_\C(X^\bullet,Y^\bullet)$ satisfying $(f_0)_\ast\delta=(f_{n+1})^\ast\rho.$

We use the same composition and identities as in ${\bf{C}}^{n+2}_\C$.
\end{itemize}

\end{definition}

\begin{definition}\cite[Definition 2.11]{HLN}
By Yoneda lemma, any extension $\delta\in\mathbb{E}(C,A)$ induces natural transformations
$$\delta_\sharp:\C(-,C)\Rightarrow\mathbb{E}(-,A)~~\text{and}~~\delta^\sharp:\C(A,-)
\Rightarrow\mathbb{E}(C,-).$$
For any $X\in\C$, these $(\delta_\sharp)_X$ and $\delta^\sharp_X$ are given as follows.
\begin{itemize}
\item[(a)] $(\delta_\sharp)_X:\C(X,C)\rightarrow\mathbb{E}(X,A)~; f\mapsto f^\ast\delta$.
\item[(b)] $\delta^\sharp_X:\C(A,X)\rightarrow\mathbb{E}(C,X)~; g\mapsto g_\ast\delta$.
\end{itemize}
We abbreviately denote $(\delta_\sharp)_X(f)$ and $\delta^\sharp_X(g)$ by
$\delta_\sharp(f)$ and $\delta^\sharp(g)$, respectively.
\end{definition}

\begin{definition}\cite[Definition 2.13]{HLN}
An \emph{$n$-exangle} is an object $\langle X^\bullet,\delta\rangle$ in ${\AE}$ that satisfies the listed conditions.
\begin{itemize}
\item[(1)]The following sequence of functors $\C^{op}\rightarrow\rm{Ab}$ is exact.
$$\C(-,X_0) \xrightarrow{~\C(-,d^X_0)~}\C(-,X_1)\xrightarrow{~\C(-,d^X_1)~} \cdots \xrightarrow{~\C(-,d^X_n)~} \C(-,X_{n+1}) \xrightarrow{~\delta_\sharp~} \mathbb{E}(-,X_0)$$
\item[(2)]The following sequence of functors $\C^{op}\rightarrow\rm{Ab}$ is exact.
$$\C(X_{n+1},-) \xrightarrow{~\C(d^X_n,-)~}\C(X_n,-)\xrightarrow{~\C(d^X_{n-1},-)~} \cdots \xrightarrow{~\C(d^X_0,-)~} \C(X_0,-) \xrightarrow{~\delta^\sharp~} \mathbb{E}(X_{n+1},-)$$
In particular, any $n$-exangle is an object in ${\AE}$. A \emph{morphism of $n$-exangles} simply means a morphism in ${\AE}$. Thus $n$-exangles forms a full subcategory of ${\AE}$.

\end{itemize}
\end{definition}

\begin{definition}\cite[Definition 2.22]{HLN}
Let $\mathfrak{s}$ be a correspondence which associates a homotopic equivalence class $\mathfrak{s}(\delta)=[_AX^\bullet_C]$ to each extension $\delta={_A}X^\bullet_C$. Such $\mathfrak{s}$ is called a \emph{realization} of $\mathbb{E}$ if it satisfies the following conditions for any $\mathfrak{s}(\delta)=[X^\bullet]$ and any $\mathfrak{s}(\rho)=[Y^\bullet]$.
\begin{itemize}
\item[($\rm{R0}$)] For any morphism of extensions $(a,c):\delta\rightarrow\rho$, there exists a morphism $f^\bullet\in{\bf{\C}}^{n+2}_{\C}(X^\bullet,Y^\bullet)$ of the form $f^\bullet=(a,f_1,\cdots,f_n,c)$. Such $f^\bullet$ is called a \emph{lift} of $(a,c)$.

\end{itemize}
 In such a case, we simply say that ``$X^\bullet$ realizes $\delta$'' whenever they satisfy $\mathfrak{s}(\delta)=[X^\bullet]$. Moreover, a realization $\mathfrak{s}$ of $\mathbb{E}$ is said to be \emph{exact} if it satisfies the following conditions.
 \begin{itemize}
\item[($\rm{R1}$)] For any $\mathfrak{s}(\delta)=[X^\bullet]$, the pair $\langle X^\bullet,\delta\rangle$ is an $n$-exangle.
\item[($\rm{R2}$)] For any object $A\in\C$, the zero element $_A0_0=0\in\mathbb{E}(0,A)$ satisfies
    $$\mathfrak{s}(_A0_0)=[A \xrightarrow{~\rm{id}_A~} A \xrightarrow{~~}0 \xrightarrow{~~} \cdots \xrightarrow{~~}0 \xrightarrow{~~} 0].$$
    Dually, $\mathfrak{s}(_00_A)=[0 \xrightarrow{~\rm~} 0 \xrightarrow{~~}0 \xrightarrow{~~} \cdots \xrightarrow{~~}A \xrightarrow{~\rm{id}_A~} A]$ holds for any $A\in\C$.

\end{itemize}
Note that the above condition $\rm{(R1)}$ does not depend on representatives of the class $[X^\bullet]$.

\end{definition}

\begin{definition}\cite[Definition 2.23]{HLN}
Let $\mathfrak{s}$ be an exact realization of $\mathbb{E}$.
\begin{itemize}
\item[($\rm1$)] An $n$-exangle $\langle X^\bullet,\delta\rangle$ is called a \emph{$\mathfrak{s}$-distinguished} $n$-exangle if it satisfies $\mathfrak{s}(\delta)=[X^\bullet]$. \\
    We often simply say \emph{distinguished $n$-exangle} when $\mathfrak{s}$ is clear from the context.
\item[($\rm2$)] An object $X^\bullet\in{\bf{C}}^{n+2}_{\C}$ is called an \emph{$\mathfrak{s}$-conflation} if it realizes some extension $\delta\in\mathbb{E}{(X_{n+1},X_0)}$.
\item[($\rm3$)] A morphism $f$ in $\C$ is called an \emph{$\mathfrak{s}$-inflation} or simply an \emph{inflation} if it admits some conflation $X^\bullet\in{\bf{C}}^{n+2}_{\C}$ satisfying $d^0_X=f$.
\item[($\rm4$)] A morphism $g$ in $\C$ is called an \emph{$\mathfrak{s}$-deflation} or simply an \emph{deflation} if it admits some conflation $X^\bullet\in{\bf{C}}^{n+2}_{\C}$ satisfying $d^n_X=g$.

\end{itemize}

\end{definition}

\begin{definition}\cite[Definition 2.27]{HLN}
For a morphism $f^\bullet\in{\bf{C}}^{n+2}_{\C}(X^\bullet,Y^\bullet)$ satisfying $f^0={\rm{id}}_A$ for some $A=X_0=Y_0$, its \emph{mapping cone} $M^\bullet_f\in{\bf{C}}^{n+2}_{\C}$ is defined to be the complex
$$X_1 \xrightarrow{~d^{M_f}_0~} X_2\oplus Y_1 \xrightarrow{~d^{M_f}_1~} X_3\oplus Y_2 \xrightarrow{~d^{M_f}_2~}\cdots \xrightarrow{~d^{M_f}_{n-1}~} X_{n+1}\oplus Y_n \xrightarrow{~d^{M_f}_n~} Y_{n+1}$$
where $d^{M_f}_0=\left[\begin{matrix}
                      -d^X_1 \\
                      f_1 \\
                    \end{matrix}\right],d^{M_f}_i=\left[\begin{matrix}
                      -d^X_{i+1}&0 \\
                      f_{i+1}&d^Y_i \\
                    \end{matrix}\right](1\leq i\leq n-1),d^{M_f}_n=\left[\begin{matrix}
                      f_{n+1}&d^Y_n \\
                    \end{matrix}\right]. $

The \emph{mapping cocone} is defined dually, for morphism $h^\bullet$ in ${\bf C}^{n+2}_{\C}$ satisfying $h_{n+1}=\rm{id}.$

\end{definition}

\begin{definition}\cite[Definition 2.32]{HLN}
An \emph{$n$-exangulated category} is a triple $(\C,\mathbb{E},\mathfrak{s})$ of additive category $\C$, additive bifunctor $\mathbb{E}:{\C}^{op}\times\C\rightarrow{\rm{Ab}}$, and its exact realization $\mathfrak{s}$, satisfying the following conditions.
\begin{itemize}
\item[($\rm{EA1}$)]Let $A \xrightarrow{~f~} B \xrightarrow{~g~} C$ be any sequence of morphism in $\C$. If both $f$ and $g$ are inflations, then so is $g\circ f$. Dually, if $f$ and $g$
are deflations, then so is $g\circ f$.
\item[($\rm{EA2}$)]For any $\rho\in\mathbb{E}(D,A)$ and $c\in\C(C,D)$, let $_A\langle X^\bullet,c^\ast\rho\rangle_C$ and $_A\langle Y^\bullet,\rho\rangle_D$ be distinduished $n$-exangles. Then $(\rm{id}_A,c)$ has a \emph{good lift} $f^\bullet$, in the sense that its mapping cone gives a distinguished $n$-exangle
    $\langle M^\bullet_f,(d^X_0)_\ast\rho\rangle$.
\item[$(\rm{{EA2})^{op}}$]Dual of ($\rm{EA2}$).
\end{itemize}
Note that the case $n=1$, a triplet $(\C, \mathbb{E},\mathfrak{s})$ is a $1$-exangulated category if and only if it is an exangulated category, see \cite[Proposition 4.3]{HLN}.
\end{definition}

\begin{example}
From \cite[Proposition 4.5]{HLN} and \cite[Proposition 4.34]{HLN}, we know that $(n+2)$-angulated categories and $n$-exact categories are $n$-exangulated categories. There are some other examples of $n$-exangulated categories which are neither $n$-exact nor $(n+2)$-angulated, see \cite[Section 6]{HLN} for more details.

\end{example}

\begin{definition}\cite[Definition 3.2]{LZ}
Let $(\C,\mathbb{E},\mathfrak{s})$ be an $n$-exangulated category.
\begin{itemize}
\item[($\rm{1}$)]An object $P\in\C$ is called \emph{projective} if for any distinguished $n$-exangle
    $$\xymatrix{
X_0 \ar[r]^-{d^X_0} & X_1 \ar[r]^-{d^X_1} & \cdots \ar[r]^-{d^X_{n-2}}&X_{n-1} \ar[r]^-{d^X_{n-1}}& X_n \ar[r]^-{d^X_{n}} & X_{n+1} \ar@{-->}[r]^{\delta}&}$$
and any morphism $c$ in $\C(P,X_{n+1})$, there exists a morphism $b\in\C(P,X_n)$ such that $d^X_n\circ b=c$. The full subcategory of projectives is denoted by $\P$. Dually, The full subcategory of injectives is denoted by $\I$.
\item[($\rm{2}$)]We say that $\C$  \emph{has enough projectives} if for any object $C\in\C$, there exists a distinguished $n$-exangle
 $$\xymatrix{
 B\ar[r]^-{\alpha_0} & P_1 \ar[r]^-{\alpha_1}& P_2 \ar[r]^-{\alpha_2} & \cdots \ar[r]^-{\alpha_{n-2}}&P_{n-1} \ar[r]^-{\alpha_{n-1}}& P_n \ar[r]^-{\alpha_{n}} & C \ar@{-->}[r]^-{\delta}&}$$
satisfying $P_1,P_2,\cdots,P_n\in\P$. Dually, we can define the notion of \emph{having enough injectives}.
\item[($\rm{3}$)]$\C$ is said to be \emph{Frobenius} if $\C$ has enough projectives and enough injectives and if moreover the projectives coincide with the injectives.

\end{itemize}

\end{definition}

\begin{remark}
\begin{itemize}
\item[($\rm{1}$)]In case $n=1$, these agree with the usual definitions
\cite[Definition 3.23, Definition 3.25 and Definition 7.1]{NP}.
\item[($\rm{2}$)] If $(\C,\mathbb{E},\mathfrak{s})$ is an $n$-exact category, then these agree with \cite[Definition 3.11, Definition 5.3 and Definition 5.5]{j}.
\item[($\rm{3}$)]If $(\C,\mathbb{E},\mathfrak{s})$ is an $n+2$-angulated category, then $\P=\I$ consists of zero objects. Moreover it always has enough projectives and enough injectives.
\end{itemize}

\end{remark}

The following are some very useful lemmas and they will be needed later on.

\begin{lemma}\textnormal{\cite[lemma 2.12]{LZ}}\label{lemma212}
Let $(\C,\mathbb{E},\mathfrak{s})$ be an $n$-exangulated category, and
$$\xymatrix{
X_0 \ar[r]^-{d^X_0} & X_1 \ar[r]^-{d^X_1} & \cdots \ar[r]^-{d^X_{n-2}}&X_{n-1} \ar[r]^-{d^X_{n-1}}& X_n \ar[r]^-{d^X_{n}} & X_{n+1} \ar@{-->}[r]^{\delta}&}$$
be a distinguished $n$-exangle. Then the following long exact sequences holds:
$$\C(-,X_0)\rightarrow\C(-,X_1)\rightarrow\cdots\rightarrow\C(-,X_{n+1})\rightarrow\mathbb{E}(-,X_0)
\rightarrow\mathbb{E}(-,X_1)\rightarrow\mathbb{E}(-,X_2);$$
$$\C(X_{n+1},-)\rightarrow\C(X_{n},-)\rightarrow\cdots\rightarrow\C(X_{0},-)\rightarrow\mathbb{E}(X_{n+1},-)
\rightarrow\mathbb{E}(X_{n},-)\rightarrow\mathbb{E}(X_{n-1},-).$$

In particular, $d^X_i$ is a weak cokernel of $d^X_{i-1}$ and $d^X_{i-1}$ is a weak kernel of  $d^X_i$, for each $i\in\{1,2,\cdots,n\}$.
\end{lemma}

Let $(\C,\mathbb{E},\mathfrak{s})$ be an $n$-exangulated category with enough projectives $\P$ and injectives $\I$. Let $X$ be any object in $\C$. It admit a distinguished $n$-exangle
$$\xymatrix{
 TX\ar[r]^-{f_0} & P_1^X \ar[r]^-{f_1}& P_2^X \ar[r]^-{f_2} & \cdots \ar[r]^-{f_{n-2}}&P_{n-1}^X \ar[r]^-{f_{n-1}}& P_n^X \ar[r]^-{f_{n}} & X \ar@{-->}[r]^-{{\delta}^X}&}$$
with $P_1^X,P_2^X,\cdots,P_n^X\in\P$ and a distinguished $n$-exangle

$$\xymatrix{
 X\ar[r]^-{g_0} & I_1^X \ar[r]^-{g_1}& I_2^X \ar[r]^-{g_2} & \cdots \ar[r]^-{g_{n-2}}&I_{n-1}^X \ar[r]^-{g_{n-1}}& I_n^X \ar[r]^-{g_{n}} & SX \ar@{-->}[r]^-{{\delta}_X}&}$$
with $I_1^X,I_2^X,\cdots,I_n^X\in\I$.

\begin{lemma}\textnormal{\cite[Lemma 3.4]{LZ}}\label{lemma213}
Let $(\C,\mathbb{E},\mathfrak{s})$ be an $n$-exangulated category. Then the following statements are equivalent for an object $P\in\C$.
\begin{itemize}
\item[$\rm{(1)}$] $\mathbb{E}(P,A)=0$ for any $A\in\C$.
\item[$\rm{(2)}$] $P$ is projective.
\item[$\rm{(3)}$] Any distinguished $n$-exangle $\xymatrix{
 A_0\ar[r]^{\alpha_0} & A_1 \ar[r]^{\alpha_1}& \cdots \ar[r]^{\alpha_{n-2}}&A_{n-1} \ar[r]^{\alpha_{n-1}}& A_n \ar[r]^{\alpha_{n}} & A_{n+1} \ar@{-->}[r]^{\delta}&}$ is split.

\end{itemize}
\end{lemma}

\begin{lemma}\textnormal{\cite[lemma 3.3]{ZW}}\label{lemma214}
Let  $(\C,\mathbb{E},\mathfrak{s})$ be an $n$-exangulated category, and

$$\xymatrix{
A_0 \ar[r]^{\alpha_0}\ar[d]^{f_0} & A_1 \ar[r]^{\alpha_1}\ar[d]^{f_1} & A_2 \ar[r]^{\alpha_2}\ar[d]^{f_2} & \cdots \ar[r]^{\alpha_{n-1}}& A_{n} \ar[r]^{\alpha_{n}}\ar[d]^{f_{n}} & A_{n+1} \ar@{-->}[r]^{\delta}\ar[d]^{f_{n+1}} &{~}\\
B_0 \ar[r]^{\beta_0} & B_1 \ar[r]^{\beta_1} & B_2 \ar[r]^{\beta_2} & \cdots \ar[r]^{\beta_{n-1}}& B_{n} \ar[r]^{\beta_{n}} & B_{n+1} \ar@{-->}[r]^{\delta'} & {~}}$$
be any morphism of distinguished $n$-exangles. Then the following are equivalent:
\begin{itemize}
\item[$\rm{(1)}$]$(f_0)_\ast\delta=f^\ast_{n+1}\delta'=0$.
\item[$\rm{(2)}$]$f_0$ factors through $\alpha_0$.
\item[$\rm{(3)}$]$f_{n+1}$ factors through $\beta_n$.

\end{itemize}
In particular, In the case $\delta=\delta'$ and $(f_0,f_1,\cdots,f_{n+1})=(1_{A_0},1_{A_1},\cdots,1_{A_{n+1}})$, $\alpha_0$ is a section if and only if $\delta=0$ if and only if $\alpha_{n}$ is a retraction\textnormal{\cite[lemma 2.15]{HLN}}.
\end{lemma}

\begin{lemma}\textnormal{\cite[Proposition 3.6]{HLN}}\label{lemma215}
Let $\langle X^\bullet,\delta\rangle$ and $\langle B^\bullet,\delta\rangle$ be distinguished $n$-exangles. Suppose that we are given a commutative square
$$\xymatrix{
 X_0\ar[r]^{d^X_0}\ar[d]^{a}&X_1\ar[d]^{b}\\
 {Y_0}\ar[r]^{d^Y_0}&Y_1 }$$
in $\C$, Then there is a morphism $f^\bullet:\langle X^\bullet,\delta\rangle\rightarrow\langle B^\bullet,\delta\rangle$ which satisfies $f_0=a$ and $f_1=b$.
\end{lemma}

\begin{lemma}\label{lemma216}
Let  $(\C,\mathbb{E},\mathfrak{s})$ be an $n$-exangulated category, and
$$\xymatrix{
A_0 \ar[r]^{\alpha_0}\ar[d]^{\varphi_0} & A_1 \ar[r]^{\alpha_1}\ar[d]^{\varphi_1} & A_2 \ar[r]^{\alpha_2}\ar[d]^{\varphi_2} & \cdots \ar[r]^{\alpha_{n-1}}& A_{n} \ar[r]^{\alpha_{n}}\ar[d]^{\varphi_{n}} & A_{n+1} \ar@{-->}[r]^{\delta}\ar[d]^{\varphi_{n+1}} &{~}\\
B_0 \ar[r]^{\beta_0} & B_1 \ar[r]^{\beta_1} & B_2 \ar[r]^{\beta_2} & \cdots \ar[r]^{\beta_{n-1}}& B_{n} \ar[r]^{\beta_{n}} & B_{n+1} \ar@{-->}[r]^{\delta'} & {~}}$$
be any morphism of distinguished $n$-exangles. Then the following are morphisms of distinguished $n$-exangles.
$$\xymatrix{A^{\bullet}\ar[d]^-{\pi^\bullet} & A_0 \ar[r]^{\alpha_0}\ar[d]^{\varphi_0} & A_1 \ar[r]^{\alpha_1}\ar[d]^{g_1} & \cdots \ar[r]^{\alpha_{n-2}}& A_{n-1} \ar[r]^{\alpha_{n-1}}\ar[d]^{g_{n-1}} & A_{n} \ar[r]^{\alpha_{n}}\ar[d]^{g_{n}}& A_{n+1} \ar@{-->}[r]^{\delta}\ar@{=}[d]^{} &{~}\\
I({\varphi^\bullet})\ar[d]^-{i^\bullet}  &B_0 \ar[r]^{e_0}\ar@{=}[d]^{} & E_1 \ar[r]^{e_1}\ar[d]^{h_1} & \cdots \ar[r]^{e_{n-2}}& E_{n-1} \ar[r]^{e_{n-1}}\ar[d]^{h_{n-1}} & E_{n} \ar[r]^{e_{n}}\ar[d]^{h_{n}}& A_{n+1} \ar@{-->}[r]^{{\varphi_0}_{\ast}\delta}\ar[d]^{\varphi_{n+1}} &{~}\\
B^{\bullet}\ar[d]^-{c^\bullet} & B_0 \ar[r]^{\beta_0}\ar[d]^{e_0} & B_1 \ar[r]^{\beta_1}\ar[d]_{\tiny{\left[\begin{matrix}
0\\
 1\\
 \end{matrix}\right]}} &\cdots \ar[r]^{\beta_{n-2}}& B_{n-1} \ar[r]^{\beta_{n-1}}\ar[d]_{\tiny{\left[\begin{matrix}
0\\
 1\\
 \end{matrix}\right]}} & B_{n} \ar[r]^{\beta_{n}}\ar[d]^{\tiny{\left[\begin{matrix}
0\\
 1\\
 \end{matrix}\right]}}& B_{n+1} \ar@{-->}[r]^{\delta'}\ar@{=}[d]^{} &{~}\\
C(\varphi^\bullet) & E_1 \ar[r]^-{\tiny{\left[\begin{matrix}
-e_1\\
 h_1\\
 \end{matrix}\right]}} & E_2\oplus B_1 \ar[r]^-{\tiny{\left[\begin{matrix}
-e_2&0\\
 h_2&\beta_1\\
 \end{matrix}\right]}} &\cdots \ar[r]^-{~}& E_{n}\oplus B_{n-1} \ar[r]^{\tiny{\left[\begin{matrix}
-e_n&0\\
 h_n&\beta_{n-1}\\
 \end{matrix}\right]}} & ~~A_{n+1}\oplus B_n\ar[r]^-{\tiny{\left[\begin{matrix}
 \varphi_{n+1}&\beta_{n}
 \end{matrix}\right]}}& ~B_{n+1} \ar@{-->}[r]^-{{e_0}_\ast\delta'} & {~}}$$
Moreover, $\underline{\varphi^\bullet}=\underline{i^\bullet\pi^\bullet}$ in $S(\C)/\R_2$.

\end{lemma}

\proof According to the assumption, we have ${\varphi_0}_{\ast}\delta={\varphi_{n+1}}^{\ast}\delta'$. Hence, By The Definition of (AE2),
 there exist two good lift morphisms $\pi^\bullet:A^\bullet\rightarrow I(\varphi^\bullet)$ and $i^\bullet:I(\varphi^\bullet)\rightarrow B^\bullet$ about distinguished $n$-exangles.  Therefore, $C(\varphi^\bullet)$ is a distinguished $n$-exangle, thus $c^\bullet:B^\bullet\rightarrow C(\varphi^\bullet)$ is morphism of distinguished $n$-exangles. It follows that the above commutative diagram holds. By Lemma \ref{lemma214}, it is clearly that $\underline{\varphi^\bullet}=\underline{i^\bullet\pi^\bullet}$ in $S(\C)/\R_2$.  \qed

\section{Identifying quotients of morphism categories as module categories }\label{3}

Our approach to understanding the category of $n$-exangles will be based on viewing them as morphism categories, which we are able to identify with certain module categories.

Let $\C$ be an additive category, For two objects $f:A\rightarrow B$ and $f':A'\rightarrow B'$ in Mor$(\C)$, we define $\cal{R}(f,f')$ (resp. ${\cal{\R}}'(f,f')$) to be the set of morphisms $(a,b)$ from $f$ to $f'$ sucu that there is some morphism $p:Y\rightarrow X'$ such that $pf'=b$ (resp. $pf=a$). Then $\cal{\R}$ and ${\cal{R}}'$ are ideals of ${\rm{Mor}}(\C)$. We denote by s-epi($\C$) (resp. s-mono($\C$)) the full subcategory of Mor($\C$) consisting of split epimorphisms (resp. split monomorphisms).

Define a functor
$\alpha:\text{Mor}(\C)\rightarrow\text{mod}(\C)$ by mapping $f:A\rightarrow B$ to $F=\text{Coker}(\C(-,f):\C(-,A)\rightarrow\C(-,B))$. The functor $\alpha$ induces the following equivalences, which are crucial in the proof of our main results.

\begin{lemma}\textnormal{\cite[Lemma 3.1,Proposition 3.3]{L}}\label{lemma31}
Let $\C$ be an additive category, then
\begin{itemize}
\item[\rm(1)] \textnormal{Mor$(\C)/\R\cong$ Mor$(\C)/$[s-epi$(\C)]\cong$ mod-$\C$.}
\item[\rm(2)] \textnormal{Mor$(\C)/\R'\cong$ Mor$(\C)/$[s-mono$(\C)]\cong$((mod-${\C}^{op})^{op}$.}
\end{itemize}
\end{lemma}

Throughout this paper, we suppose that $(\C,\mathbb{E},\mathfrak{s})$ is an $n$-exangulated category. Assume that $\M$ is a full subcategory of $\C$. For convenience, we fix some notations. We denote by $\mathfrak{s}$-def$(\M)$
(resp. $\mathfrak{s}$-inf$(\M)$) the full subcategory of Mor($\C$) consisting of $\mathfrak{s}$-deflations (resp. $\mathfrak{s}$-inflations). Recall that the full subcategory of $\mathfrak{s}$-def$(\M)$ consisting of split epimorphisms (resp. split monomorphisms) is denoted by s-epi($\M$) (resp. s-mono($\M$)). We denote by sp-epi($\M$)(resp. si-mono$(\M)$) the full subcategory of $\mathfrak{s}$-def$(\M)$ consisting of $(\xymatrix{M \ar[r]^-{1} & M})\oplus(P\rightarrow M')$ (resp. ($\xymatrix{M \ar[r]^-{1} & M})\oplus(M'\oplus I)$) with $P\in\P$ (resp. $I\in\I$).

\begin{lemma}\label{lemma32}
Let $(\C,\mathbb{E},\mathfrak{s})$ be an $n$-exangulated category with enough projectives $\P$
and $\M$ be a full subcategory of $\C$ containing $\P$. Suppose that the following
$$\xymatrix{
A_0 \ar[r]^{\alpha_0}\ar[d]^{g_0} & A_1 \ar[r]^{\alpha_1}\ar[d]^{g_1} & A_2 \ar[r]^{\alpha_2}\ar[d]^{g_2} & \cdots \ar[r]^{\alpha_{n-2}}& A_{n-1} \ar[r]^{\alpha_{n-1}}\ar[d]^{g_{n-1}} & M_1 \ar[r]^{f}\ar[d]^{a} & M_2 \ar@{-->}[r]^{\delta}\ar[d]^{b} &{~}\\
A_0' \ar[r]^{\alpha_0'} & A_1' \ar[r]^{\alpha_1'} & A_2' \ar[r]^{\alpha_2'} & \cdots \ar[r]^{\alpha_{n-2}'}& A_{n-1}' \ar[r]^{\alpha_{n-1}'} & M_1' \ar[r]^{f'} & M_2' \ar@{-->}[r]^{\delta'} & {~}}$$
is a morphism of distinguished $n$-exangles. Then
\begin{itemize}
\item[$\rm{(1)}$] The following statements are equivalent.
\subitem{$\rm{(a)}$} The morphism $\underline{b}$ factors through $\underline{f'}$ in $\M/[\P]$.
\subitem{$\rm{(b)}$} The morphism $b$ factors through $f'$.
\subitem{$\rm{(c)}$}The morphism $(a,b)$ factors through some object in \textnormal{s-epi}$(\M)$.
\item[$\rm{(2)}$] The following statements are equivalent.
\subitem{$\rm{(a)}$} The morphism $\underline{a}$ factors through $\underline{f}$ in $\M/[\P]$.
\subitem{$\rm{(b)}$}The morphism $(a,b)$ factors through some object in \textnormal{sp-epi}$(\M)$.
\end{itemize}
\proof (1) By Lemma \ref{lemma31}, $\rm(b)$ and $\rm(c)$ are equivalent. It is clearly that $\rm(b)$ implies $\rm(a)$, we only need to show that $\rm(a)$ implies $\rm(b)$. This is an adaptation of the proof of {\cite[Lemma 2.11]{Li}}. Assume that there is a morphism $\underline{p}:M_2\rightarrow M_1'$ sucu that $\underline{b}=\underline{f'p}$. Hence there exists an object $P\in\P$ and two morphisms $u:M_2\rightarrow P$ and $v:P\rightarrow M_2'$ such that $b-f'p=vu$. Since $f'$ is an $\mathfrak{s}$-deflation, there exists a morphism
$w:P\rightarrow M_1'$ such that $v=f'w$. Therefore, $b=f'(p+wu)$.

(2) The proof is similar to {\cite[Lemma 3.1(2)]{Li}}, we omit it. \qed

\end{lemma}

\begin{lemma}\label{le33}
Let  $(\C,\mathbb{E},\mathfrak{s})$ be an $n$-exangulated category. Assume that
 $$\xymatrix{
 A\ar[r]^{\beta_0} & B_1 \ar[r]^{\beta_1}& B_2 \ar[r]^{\beta_2} & \cdots \ar[r]^{\beta_{n-2}}&B_{n-1} \ar[r]^{\beta_{n-1}}& B_n \ar[r]^{\beta_{n}} & D \ar@{-->}[r]^{\theta}&}$$
 is a distinguished $n$-exangle, $h:C\rightarrow D$ is a morphism and
 $$\xymatrix{
 A\ar[r]^{\alpha_0} & A_1 \ar[r]^{\alpha_1}& A_2 \ar[r]^{\alpha_2} & \cdots \ar[r]^{\alpha_{n-2}}&A_{n-1} \ar[r]^{\alpha_{n-1}}& A_n \ar[r]^{\alpha_{n}} & C \ar@{-->}[r]^{h^\ast\theta}&}$$ is a distinguished $n$-exangle, then there is a morphism $f^\bullet:=(1_A,f_1,\cdots,f_n,h)$ of distinguished $n$-exangles
 $$\xymatrix{
A \ar[r]^{\alpha_0}\ar@{=}[d]& A_1 \ar[r]^{\alpha_1}\ar[d]^{f_1} & A_2 \ar[r]^{\alpha_2}\ar[d]^{f_2} & \cdots \ar[r]^{\alpha_{n-1}}& A_{n} \ar[r]^{\alpha_{n}}\ar[d]^{f_{n}} & C \ar@{-->}[r]^{h^\ast\theta}\ar[d]^{h} &{~}\\
A\ar[r]^{\beta_0} & B_1 \ar[r]^{\beta_1} & B_2 \ar[r]^{\beta_2} & \cdots \ar[r]^{\beta_{n-1}}& B_{n} \ar[r]^{\beta_{n}} & D \ar@{-->}[r]^{\theta} & {~}}$$
and moreover, $\xymatrix{
 A_1\ar[r]^-{\gamma_0} & A_2\oplus B_1 \ar[r]^-{\gamma_1}& A_3\oplus B_2 \ar[r]^-{\gamma_2} & \cdots \ar[r]^-{\gamma_{n-2}}& C\oplus B_n \ar[r]^-{\gamma_{n}} & D \ar@{-->}[r]^{{\alpha_0}_{\ast}\theta}&}$ is a distinguished $n$-exangle

where $\gamma_0=\left[\begin{matrix}
                      -\alpha_1 \\
                      f_1 \\
                    \end{matrix}\right],\gamma_i=\left[\begin{matrix}
                      -\alpha_{i+1}&0 \\
                      f_{i+1}&\beta_i \\
                    \end{matrix}\right](1\leq i\leq n-1),\gamma_n=\left[\begin{matrix}
                      h&\beta_n \\
                    \end{matrix}\right]. $
\proof This follows from the definition of $\rm{(EA2)}$. \qed                    
\end{lemma}

\begin{theorem}\label{th34}
Let $(\C,\mathbb{E},\mathfrak{s})$ be an $n$-exangulated category and $\M$ be a full subcategory of $\C$.
\begin{itemize}
\item[$\rm{(i)}$] If $\C$ has enough projectives $\P$ and $\M$ contains $\P$, then we have the following equivalences:
$$\mathfrak{s}\textnormal{-def}(\M)/[\textnormal{s-epi}(\M)]\cong\textnormal{mod-}(\M/[\P]),$$ $$\mathfrak{s}\textnormal{-def}(\M)/[\textnormal{sp-epi}(\M)]\cong
{(\textnormal{mod-}{{(\M/[\P])}^{op}})}^{op}.$$
\item[$\rm{(ii)}$] If $\C$ has enough injectives $\I$ and $\M$ contains $\I$, then we have the following equivalences:
$$\mathfrak{s}\textnormal{-inf}(\M)/[\textnormal{s-mono}(\M)]\cong
{(\textnormal{mod-}{{(\M/[\I])}^{op}})}^{op},$$ $$\mathfrak{s}\textnormal{-inf}(\M)/[\textnormal{si-mono}(\M)]\cong\textnormal{mod-}(\M/[\I]).$$
\end{itemize}
\proof  We only prove the $(\rm i)$ since the statement $(\rm ii)$ can be proved dually.

 By Lemma \ref{lemma32}, we naturally give an additive functor
$$F:\mathfrak{s}\textnormal{-def}(\M)\rightarrow \textnormal{Mor}(\M/[\P]), (\xymatrix{M_1 \ar[r]^-{f} & M_2})\mapsto\xymatrix{M_1 \ar[r]^-{\underline{f}} & M_2}.$$
For any object $\underline{f}:M_1\rightarrow M_2$ in Mor$(\M/[\P])$, Since $\C$ has enough projectives, there is an $\mathfrak{s}$-deflation $\pi:P\rightarrow M_2$ with $P\in\P$. It follows from Lemma \ref{le33} that $(f,\pi):M_1\oplus P\rightarrow M_2$ is an $\mathfrak{s}$-deflation and $\underline{(f,\pi)}\cong\underline{f}$. It implies that $F(f,\pi)\cong\underline{f}$. Hence $F$ is dense.

Suppose that $f:M_1\rightarrow M_2$ and $f':M_1'\rightarrow M_2'$ are objects in $\mathfrak{s}$-def$(\M)$ and $(\underline{a},\underline{b})$ is a morphism
in $\M/[P]$ from $\underline{f}$ to $\underline{f'}$. Thus, we have $\underline{bf}=\underline{f'a}$. Hence there exists an object $P\in\P$ and two morphisms $u:M_1\rightarrow P$ and $v:P\rightarrow M_2'$ such that $bf-f'a=vu$. Since $f'$ is an $\mathfrak{s}$-deflation, there exists a morphism $w:P\rightarrow M_1'$ such that $v=f'w$. Therefore, $bf=f'(a+wu)$. This shows that $F(a+wu,b)=(\underline{a},\underline{b})$. Hence $F$ is full.

The functor $F$ induces a full and dense functor $\tilde{F}:\mathfrak{s}\textnormal{-def}(\M)\rightarrow\textnormal{Mor}(\M/[\P])/{\R}$. According to Lemma \ref{lemma32}(1), we have $\mathfrak{s}\textnormal{-def}(\M)/[\textnormal{s-epi}(\M)]\cong\textnormal{Mod-}(\M/[\P])/\R$. Therefore, by Lemma \ref{lemma31}(1) shows that $\mathfrak{s}\textnormal{-def}(\M)/[\textnormal{s-epi}(\M)]\cong\textnormal{mod-}(\M/[\P])$.

The functor $F$ induces a full and dense functor $\hat{F}:\mathfrak{s}\textnormal{-def}(\M)\rightarrow\textnormal{Mor}(\M/[\P])/{\R}'$. According to Lemma \ref{lemma32}(2), we have $\mathfrak{s}\textnormal{-def}(\M)/[\textnormal{sp-epi}(\M)]\cong\textnormal{Mod-}(\M/[\P])/\R'$. Therefore, by Lemma \ref{lemma31}(2) shows that $\mathfrak{s}\textnormal{-def}(\M)/[\textnormal{sp-epi}(\M)]\cong
(\textnormal{mod-}(\M/[\P])^{op})^{op}$. This concludes the proof  \qed

\end{theorem}

\section{Abelian quotients of categories of $n$-exangles}\label{4}

In this section, we assume that $(\C,\mathbb{E},\mathfrak{s})$ is an $n$-exangulated category. We denote by $S(\C)$ the category of $n$-exangles, where the objects are the distinguished $n$-exangles in $\C$

$$A^\bullet=(\xymatrix{
 A_0\ar[r]^{\alpha_0} & A_1 \ar[r]^{\alpha_1}& A_2 \ar[r]^{\alpha_2} & \cdots \ar[r]^{\alpha_{n-2}}&A_{n-1} \ar[r]^{\alpha_{n-1}}& A_n \ar[r]^{\alpha_{n}} & A_{n+1} \ar@{-->}[r]^{\delta}&})$$
and the morphisms from $A^\bullet$ to $B^\bullet$ is defined to be the $\varphi\bullet=(\varphi_0,\varphi_1,\cdots,\varphi_{n+1})$ such that the following diagram is commutative
$$\xymatrix{
A_0 \ar[r]^{\alpha_0}\ar[d]^{\varphi_0} & A_1 \ar[r]^{\alpha_1}\ar[d]^{\varphi_1} & A_2 \ar[r]^{\alpha_2}\ar[d]^{\varphi_2} & \cdots \ar[r]^{\alpha_{n-1}}& A_{n} \ar[r]^{\alpha_{n}}\ar[d]^{\varphi_{n}} & A_{n+1} \ar@{-->}[r]^{\delta}\ar[d]^{\varphi_{n+1}} &{~}\\
B_0 \ar[r]^{\beta_0} & B_1 \ar[r]^{\beta_1} & B_2 \ar[r]^{\beta_2} & \cdots \ar[r]^{\beta_{n-1}}& B_{n} \ar[r]^{\beta_{n}} & B_{n+1} \ar@{-->}[r]^{\delta'} & {~}}$$
and ${\varphi_0}_{\ast}\delta={\varphi_{n+1}}^{\ast}\delta'$. Let $A^\bullet$ and $B^\bullet$ be two distinguished $n$-exangles, we denote by $\R_2(A^\bullet,B^\bullet)$ (resp. $\R_1'(A^\bullet,B^\bullet)$) the class of morphism $\varphi^\bullet:A^\bullet\rightarrow B^\bullet$ such that $\varphi_{n+1}$ factors through $\beta_{n}$ (resp. $\varphi_{0}$ factors through $\alpha_{0}$). It is obvious that  $\R_2$ and $\R_1'$ are ideals of distinguished $n$-exangles in $\C$. Note that by Lemma \ref{lemma214}, we have $\R_2=\R_1'$.

In this section, we will realize quotients of categories of $n$-exangles as module categories and we will prove that for general $n$-exangulated category $(\C,\mathbb{E},\mathfrak{s})$, the quotient category $S(\C)/\R_2$ always has an abelian structure.

\begin{theorem}\label{th41}
Let $(\C,\mathbb{E},\mathfrak{s})$ be an $n$-exangulated category.
\begin{itemize}
\item[$\rm{(i)}$] If $\C$ has enough projectives, then we have the following equivalence:
$$S(\C)/\R_2\cong\textnormal{mod-}(\C/[\P])$$

\item[$\rm{(ii)}$] If $\C$ has enough injectives $\I$, then we have the following equivalence:
$$S(\C)/\R_2\cong{(\textnormal{mod-}{{(\C/[\I])}^{op}})}^{op}$$
\end{itemize}
\proof (1) According to Lemma \ref{lemma32}, we have $S(\C)/\R_2\cong\mathfrak{s}\textnormal{-def}(\C)/[\textnormal{s-epi}(\C)].$ Consider the case of $\M=\C$ in Theorem \ref{th34}(1), we can directly conclude that $S(\C)/\R_2\cong\textnormal{mod-}(\C/[\P])$.

(2) Due to $\R_2=\R_1'$, Consider the case of $\M=\C$ in Theorem \ref{th34}(2)
we have $S(\C)/\R_2\cong S(\C)/\R_1'\cong\mathfrak{s}\textnormal{-inf}(\M)/[\textnormal{s-mono}(\M)]
\cong{(\textnormal{mod-}{{(\C/[\I])}^{op}})}^{op}$.  \qed
\end{theorem}

\begin{lemma}\label{lemma42}
Let $(\C,\mathbb{E},\mathfrak{s})$ be an $n$-exangulated category. Assume that the following
$$\xymatrix{
A^{\bullet}\ar[d]^-{\varphi^\bullet} & A_0 \ar[r]^{\alpha_0}\ar[d]^{\varphi_0} & A_1 \ar[r]^{\alpha_1}\ar[d]^{\varphi_1} & A_2 \ar[r]^{\alpha_2}\ar[d]^{\varphi_2} & \cdots \ar[r]^{\alpha_{n-1}}& A_{n} \ar[r]^{\alpha_{n}}\ar[d]^{\varphi_{n}} & A_{n+1} \ar@{-->}[r]^{\delta}\ar[d]^{\varphi_{n+1}} &{~}\\
B^{\bullet} & B_0 \ar[r]^{\beta_0} & B_1 \ar[r]^{\beta_1} & B_2 \ar[r]^{\beta_2} & \cdots \ar[r]^{\beta_{n-1}}& B_{n} \ar[r]^{\beta_{n}} & B_{n+1} \ar@{-->}[r]^{\delta'} & {~}}$$
is a morphism of distinguished $n$-exangles. Then
\begin{itemize}
\item[$\rm{(1)}$] The following statements are equivalent.
\subitem{$\rm{(a)}$} $\underline{\varphi_\bullet}=0$ in $S(\C)/\R_2$.
\subitem{$\rm{(b)}$} $\varphi_0$ factors through $\alpha_0$.
\subitem{$\rm{(c)}$} $\varphi_{n+1}$ factors through $\beta_n$.

\item[$\rm{(2)}$] When taking $n$ as an even number, the following are morphisms of distinguished $n$-exangles.
$$\xymatrix{K(\varphi^{\bullet})\ar[d]^-{k^\bullet} & A_0 \ar[r]^-{\tiny{\left[\begin{matrix}
-\alpha_0\\
 \varphi_0\\
 \end{matrix}\right]}}\ar@{=}[d]^{} & A_1\oplus B_0 \ar[r]^-{\tiny{\left[\begin{matrix}
-\alpha_1&0\\
 g_2&e_0\\
 \end{matrix}\right]}}\ar[d]^{\tiny{\left[\begin{matrix}
-1&0\\
 \end{matrix}\right]}} & \cdots \ar[r]^{~}& A_{n-1}\oplus E_{n-2} \ar[r]^-{\tiny{\left[\begin{matrix}
-\alpha_{n-1}&0\\
 g_{n-1}&e_{n-2}\\
 \end{matrix}\right]}}\ar[d]^{\tiny{\left[\begin{matrix}
(-1)^{n-1}&0\\
 \end{matrix}\right]}} & A_{n}\oplus E_{n-1} \ar[r]^-{\tiny{\left[\begin{matrix}
g_n&e_{n-1}\\
 \end{matrix}\right]}}\ar[d]^{\tiny{\left[\begin{matrix}
(-1)^{n}&0\\
 \end{matrix}\right]}}& E_n \ar@{-->}[r]^{{e_n}^{\ast}\delta}\ar[d]^{(-1)^{n}e_n} &{~}\\
A^{\bullet}\ar[d]^-{\pi^\bullet} & A_0 \ar[r]^{\alpha_0}\ar[d]^{\varphi_0} & A_1 \ar[r]^{\alpha_1}\ar[d]^{g_1} & \cdots \ar[r]^{\alpha_{n-2}}& A_{n-1} \ar[r]^{\alpha_{n-1}}\ar[d]^{g_{n-1}} & A_{n} \ar[r]^{\alpha_{n}}\ar[d]^{g_{n}}& A_{n+1} \ar@{-->}[r]^{\delta}\ar@{=}[d]^{} &{~}\\
I({\varphi^\bullet})\ar[d]^-{i^\bullet}  &B_0 \ar[r]^{e_0}\ar@{=}[d]^{} & E_1 \ar[r]^{e_1}\ar[d]^{h_1} & \cdots \ar[r]^{e_{n-2}}& E_{n-1} \ar[r]^{e_{n-1}}\ar[d]^{h_{n-1}} & E_{n} \ar[r]^{e_{n}}\ar[d]^{h_{n}}& A_{n+1} \ar@{-->}[r]^{{\varphi_0}_{\ast}\delta}\ar[d]^{\varphi_{n+1}} &{(4.1)}\\
B^{\bullet}\ar[d]^-{c^\bullet} & B_0 \ar[r]^{\beta_0}\ar[d]^{e_0} & B_1 \ar[r]^{\beta_1}\ar[d]^{\tiny{\left[\begin{matrix}
0\\
 1\\
 \end{matrix}\right]}} &\cdots \ar[r]^{\beta_{n-2}}& B_{n-1} \ar[r]^{\beta_{n-1}}\ar[d]^{\tiny{\left[\begin{matrix}
0\\
 1\\
 \end{matrix}\right]}} & B_{n} \ar[r]^{\beta_{n}}\ar[d]^{\tiny{\left[\begin{matrix}
0\\
 1\\
 \end{matrix}\right]}}& B_{n+1} \ar@{-->}[r]^{\delta'}\ar@{=}[d]^{} &{~}\\
C(\varphi^\bullet) & E_1 \ar[r]^-{\tiny{\left[\begin{matrix}
-e_1\\
 h_1\\
 \end{matrix}\right]}} & E_2\oplus B_1 \ar[r]^-{\tiny{\left[\begin{matrix}
-e_2&0\\
 h_2&\beta_1\\
 \end{matrix}\right]}} &\cdots \ar[r]^-{~}& E_{n}\oplus B_{n-1} \ar[r]^{\tiny{\left[\begin{matrix}
-e_n&0\\
 h_n&\beta_{n-1}\\
 \end{matrix}\right]}} & A_{n+1}\oplus B_n \ar[r]^-{\tiny{\left[\begin{matrix}
 \varphi_{n+1}&\beta_{n}
 \end{matrix}\right]}}& B_{n+1} \ar@{-->}[r]^{{e_0}_\ast\delta'} & {~}}$$
 Moreover, $\underline{\varphi^\bullet}=\underline{i^\bullet\pi^\bullet}$ in $S(\C)/\R_2$.

\item[$\rm{(3)}$] The following statements are equivalent.
\subitem{$\rm{(a)}$} $\underline{\varphi^\bullet}$ is epimorphism in $S(\C)/\R_2$.
\subitem{$\rm{(b)}$} ${\left[\begin{matrix}
 \varphi_{n+1}&\beta_{n}
 \end{matrix}\right]}:A_{n+1}\oplus B_n\rightarrow B_{n+1}$ is a retraction.
\end{itemize}
\proof (1) It can be immediately obtained from Lemma \ref{lemma213}.

   (2) The duall of Lemma \ref{le33} imply that $K(\varphi^\bullet)$ is a distinguished $n$-exangle, therefore, $k^\bullet:K(\varphi^\bullet)\rightarrow A^\bullet$ is a morphism of distinguished $n$-exangle. Moreover, by Lemma \ref{lemma216}, the above diagram is commutative and $\underline{\varphi^\bullet}=\underline{i^\bullet\pi^\bullet}$ in $S(\C)/\R_2$.

  (3) $\rm(a)$ implies $\rm(b)$. Since there is a morphism $q={\left[\begin{matrix}
 1\\
 0\\
 \end{matrix}\right]}:A_{n+1}\rightarrow A_{n+1}\oplus B_n$ such that $\varphi_{n+1}=\left[\begin{matrix}
 \varphi_{n+1}&\beta_{n}
 \end{matrix}\right]q$, hence, by Lemma \ref{lemma216}, $\underline{c^\bullet\varphi^\bullet}=\underline{c^\bullet i^\bullet\pi^\bullet}=0$. By hypothesis, $\underline{\varphi^\bullet}$ is a epimorphism, therefore, $\underline{c^\bullet}=0$. It follows from  Lemma \ref{lemma42}(1) that ${\left[\begin{matrix}
 \varphi_{n+1}&\beta_{n}
 \end{matrix}\right]}:A_{n+1}\oplus B_n\rightarrow B_{n+1}$ is a retraction.

 $\rm(b)$ implies $\rm(a)$. Suppose that there exists a morphism $\left[\begin{matrix}
 \varphi_{n+1}'\\
 \beta_{n}'\\
 \end{matrix}\right]: B_{n+1}\rightarrow A_{n+1}\oplus B_n$ such that $\left[\begin{matrix}
 \varphi_{n+1}& \beta_{n}\\
 \end{matrix}\right]\left[\begin{matrix}
 \varphi_{n+1}'\\
 \beta_{n}'\\
 \end{matrix}\right]=\varphi_{n+1}\varphi_{n+1}'+\beta_n\beta_n'=1$. For any morphism $\psi^\bullet: C^\bullet\rightarrow A^\bullet$ of distinguished $n$-exangles such that $\underline{\psi^\bullet\varphi^\bullet}=0$, we need to show that $\underline{\psi^\bullet}=0$.
 $$\xymatrix{
A^{\bullet}\ar[d]^-{\varphi^\bullet} & A_0 \ar[r]^{\alpha_0}\ar[d]^{\varphi_0} & A_1 \ar[r]^{\alpha_1}\ar[d]^{\varphi_1} & A_2 \ar[r]^{\alpha_2}\ar[d]^{\varphi_2} & \cdots \ar[r]^{\alpha_{n-1}}& A_{n} \ar[r]^{\alpha_{n}}\ar[d]^{\varphi_{n}} & A_{n+1} \ar@{-->}[r]^{\delta}\ar[d]^{\varphi_{n+1}} &{~}\\
B^{\bullet}\ar[d]^-{\psi^\bullet} & B_0 \ar[r]^{\beta_0}\ar[d]^{\psi_0} & B_1 \ar[r]^{\beta_1}\ar[d]^{\psi_1} & B_2 \ar[r]^{\beta_2}\ar[d]^{\psi_2} & \cdots \ar[r]^{\beta_{n-1}}& B_{n} \ar[r]^{\beta_{n}}\ar[d]^{\psi_{n}} & B_{n+1} \ar@{-->}[r]^{\delta'}\ar[d]^{\psi_{n+1}} &{~}\\
C^\bullet &C_0 \ar[r]^{\gamma_0} & C_1 \ar[r]^{\gamma_1} & C_2 \ar[r]^{\gamma_2} & \cdots \ar[r]^{\gamma_{n-1}}& C_{n} \ar[r]^{\gamma_{n}} & C_{n+1} \ar@{-->}[r]^{\delta''} & {~}}$$

 Indeed, by Lemma \ref{lemma42}(1), there exists a morphism $h_n:A_{n+1}\rightarrow C_n$ such that $\psi_{n+1}\varphi_{n+1}=\gamma_nh_n$. Therefore, there exists a morphism $u=\left[\begin{matrix}
 h_{n}& \psi_{n}
 \end{matrix}\right]\left[\begin{matrix}
 \varphi_{n+1}'\\
 \beta_{n}'\\
 \end{matrix}\right]:B_{n+1}\rightarrow C_n$ such that
 $$\gamma_nu=\gamma_n\left[\begin{matrix}
 h_{n}& \psi_{n}
 \end{matrix}\right]\left[\begin{matrix}
 \varphi_{n+1}'\\
 \beta_{n}'\\
 \end{matrix}\right]=\gamma_n(h_n\psi_{n+1}'+\psi_n\beta_n')=
 \psi_{n+1}(\varphi_{n+1}\varphi_{n+1}'+\beta_n\beta_n')=\psi_{n+1}.$$
 This shows that $\underline{\psi^\bullet}=0$. \qed
\end{lemma}

\begin{remark}\label{re43}
Let $\C$ be $n$-exangulated category. we need to require $n$ to even in order for The dual conclusion of Lemma \ref{lemma42}(3) to holds. However, when $\C$ is $n$-exact category, $n$ does not to be odd or even, the conclusion Lemma \ref{lemma42}(3) and the duality of Lemma \ref{lemma42}(3) still holds.

\end{remark}

If $\C$ has enough projectives $\P$, then we have the following equivalence:
$S(\C)/\R_2\cong\textnormal{mod-}(\C/[\P])$ by Theorem \ref{th41}. The following result implies that the quotient category $S(\C)/\R_2$ always has an abelian structure when taking $n$ as an even number.

\begin{theorem}\label{th44}
Let $(\C,\mathbb{E},\mathfrak{s})$ be an $n$-exangulated category. When taking $n$ as an even number, the quotient category $S(\C)/\R_2$ is an abelian category.
\proof  Suppose that $\varphi^\bullet:A^\bullet\rightarrow B^\bullet$ is a morphism in $S(\C)$. As notations in daigram $\rm(4.1)$, we claim that $\underline{k^\bullet}:K(\varphi^\bullet)\rightarrow A^\bullet$ is a kernel of $\varphi^\bullet.$ By the dual of Lemma \ref{lemma42}(3), $\underline{k^\bullet}$ is a monomorphism.

Since $\varphi_0=\left[\begin{matrix}
 0&1
 \end{matrix}\right]\left[\begin{matrix}
 -\alpha_0\\
\varphi_0\\
 \end{matrix}\right]$, by Lemma \ref{lemma42}(1), $\underline{\varphi^\bullet k^\bullet}=0$. Suppose that there is a morphism $\psi^\bullet:C^\bullet\rightarrow A^\bullet$ such that $\underline{\varphi^\bullet\psi^\bullet}=0$.
  $$\xymatrix{
C^{\bullet}\ar[d]^-{\psi^\bullet} & C_0 \ar[r]^{\gamma_0}\ar[d]^{\psi_0} & C_1 \ar[r]^{\gamma_1}\ar[d]^{\psi_1} & C_2 \ar[r]^{\gamma_2}\ar[d]^{\psi_2} & \cdots \ar[r]^{\gamma_{n-1}}& C_{n} \ar[r]^{\gamma_{n}}\ar[d]^{\psi_{n}} & C_{n+1} \ar@{-->}[r]^{\delta''}\ar[d]^{\psi_{n+1}} &{~}\\
A^{\bullet}\ar[d]^-{\varphi^\bullet} & A_0 \ar[r]^{\alpha_0}\ar[d]^{\varphi_0} & A_1 \ar[r]^{\alpha_1}\ar[d]^{\varphi_1} & A_2 \ar[r]^{\alpha_2}\ar[d]^{\varphi_2} & \cdots \ar[r]^{\alpha_{n-1}}& A_{n} \ar[r]^{\alpha_{n}}\ar[d]^{\varphi_{n}} & A_{n+1} \ar@{-->}[r]^{\delta}\ar[d]^{\varphi_{n+1}} &{~}\\
B^\bullet &B_0 \ar[r]^{\beta_0} & B_1 \ar[r]^{\beta_1} & B_2 \ar[r]^{\beta_2} & \cdots \ar[r]^{\beta_{n-1}}& B_{n} \ar[r]^{\beta_{n}} & B_{n+1} \ar@{-->}[r]^{\delta'} & {~}}$$
By Lemma \ref{lemma42}(1), there exists a morphism $u_1:C_1\rightarrow B_0$ such that $\varphi_0\psi_0=u_1\gamma_0$. Hence, there is a morphism $\left[\begin{matrix}
 -\psi_1\\
 u_1\\
 \end{matrix}\right]:C_1\rightarrow A_1\oplus B_0$ such that $\left[\begin{matrix}
 -\psi_1\\
 u_1\\
 \end{matrix}\right]\gamma_0=\left[\begin{matrix}
 -\alpha_0\\
 \varphi_0\\
 \end{matrix}\right]\psi_0$, Therefore, by Lemma \ref{lemma215}, we have the following commutative diagram:
 $$\xymatrix{
C^{\bullet}\ar[d]^-{} & C_0 \ar[r]^{\gamma_0}\ar[d]^{\psi_0} & C_1 \ar[r]^{\gamma_1}\ar[d]^{\tiny{\left[\begin{matrix}
 -\psi_1\\
 u_1\\
 \end{matrix}\right]}} & C_2 \ar[r]^{\gamma_2}\ar@{-->}[d]^{\theta_2} & \cdots \ar[r]^{\gamma_{n-1}}& C_{n} \ar[r]^{\gamma_{n}}\ar@{-->}[d]^{\theta_{n}} & C_{n+1} \ar@{-->}[r]^{\delta''}\ar@{-->}[d]^{\theta_{n+1}} &{~}\\
K(\varphi^\bullet) &A_0 \ar[r]^-{\tiny{\left[\begin{matrix}
 -\alpha_0\\
 \varphi_0\\
 \end{matrix}\right]}} & A_1\oplus B_0 \ar[r]^{\beta_1} & A_2\oplus E_1 \ar[r]^{\beta_2} & \cdots \ar[r]^{\beta_{n-1}}& A_{n}\oplus E_{n-1} \ar[r]^{\beta_{n}} & E_{n} \ar@{-->}[r]^{{e_n}^\ast\delta} & {~}}$$
 It is clearly that $\underline{\psi^\bullet}=\underline{k^\bullet\theta^\bullet}$.

Dually, we can show that $\underline{c^\bullet}:B^\bullet\rightarrow C(\varphi^\bullet)$ is a cokernel of $\underline{\varphi^\bullet}$.

We need to show that Coker(Ker$(\underline{\varphi^\bullet}))\cong\textnormal{Ker(Cokernel}
(\underline{\varphi^\bullet}))$, that is, $\textnormal{Coker}(\underline{k^\bullet})\cong\textnormal{Ker}(\underline{c^\bullet})$.
Based on the above proof, we know that
$$\textnormal{Coker}(\underline{k^{\bullet}}):=\xymatrix@C=2cm@R=2cm
{ A_{1}\oplus B_0 \ar[r]^{\tiny{\left[\begin{matrix}
0&1&0\\
 0&0&0\\
0&0&e_{n-1}\\
 \end{matrix}\right]}} & ~A_1\oplus A_2\oplus E_1 \ar[r]^-{\tiny{\left[\begin{matrix}
0&1&0\\
 0&0&0\\
0&0&e_{n-1}\\
 \end{matrix}\right]}}& A_2\oplus A_3\oplus E_2~~~\ar[r]^{\tiny{\left[\begin{matrix}
0&1&0\\
 0&0&0\\
0&0&e_{n-1}\\
 \end{matrix}\right]}} &\cdots }$$
 $$\xymatrix@C=2cm@R=2cm{\cdots\ar[r]^-{\tiny{\left[\begin{matrix}
0&1&0\\
 0&0&0\\
0&0&e_{n-2}\\
 \end{matrix}\right]}} & A_{n-1}\oplus A_n\oplus E_{n-1} \ar[r]^{\small{\left[\begin{matrix}
0&1&0\\
 0&0&0\\
0&0&e_{n-1}\\
 \end{matrix}\right]}} & A_n\oplus E_n \ar[r]^{\tiny{\left[\begin{matrix}
0&e_{n}\\
 \end{matrix}\right]}} & A_{n+1} \ar@{-->}[r]^{\tiny{\left[\begin{matrix}
 -\alpha_0\\
 \varphi_0\\
 \end{matrix}\right]}_{\ast}\delta} & {~}}.$$
Then we have the following commutative diagram from $\textnormal{Coker}(\underline{k^{\bullet}})$ to $I(\varphi^\bullet)$.
 $$\xymatrix@C=2cm@R=1.8cm{
A_1\oplus B_0 \ar[r]^-{\tiny{\left[\begin{matrix}
 -1&0\\
 \alpha_1&0\\
 -g_1&-e_0\\
 \end{matrix}\right]}} \ar@{=}[d] & A_1\oplus A_2\oplus E_1 \ar[r]^-{\tiny{\left[\begin{matrix}
 \alpha_1&1&0\\
 0&\alpha_2&0\\
 0&-g_2&-e_1\\
 \end{matrix}\right]}}\ar[d]^{\tiny{\left[\begin{matrix}
 -1&0&0\\
 \alpha_1&1&0\\
 g_1&0&-1\\
 \end{matrix}\right]}} & A_2\oplus B_3\oplus E_2\ar[r]^-{\tiny{\left[\begin{matrix}
 \alpha_2&-1&0\\
 0&\alpha_3&0\\
 0&-g_3&-e_2\\
 \end{matrix}\right]}}\ar@{-->}[d]^{\tiny{\left[\begin{matrix}
 1&0&0\\
 \alpha_2&-1&0\\
 g_2&0&1\\
 \end{matrix}\right]}} & \cdots\\
 A_1\oplus B_0 \ar[r]^-{\tiny{\left[\begin{matrix}
 1&0\\
 0&0\\
 0&e_0\\
 \end{matrix}\right]}} &  A_1\oplus A_2\oplus E_1  \ar[r]^-{\tiny{\left[\begin{matrix}
 0&1&0\\
 0&0&0\\
 0&0&e_1\\
 \end{matrix}\right]}} &A_2\oplus B_3\oplus E_2 \ar[r]^-{\tiny{\left[\begin{matrix}
 0&1&0\\
 0&0&0\\
 0&0&e_2\\
 \end{matrix}\right]}}& \cdots }$$
 $$\xymatrix@C=2cm@R=1.8cm{
 \cdots{\ar[r]^-{\tiny{\left[\begin{matrix}
 \alpha_{n-2}&(-1)^{n-1}&0\\
 0&\alpha_{n-1}&0\\
 0&-g_{n-1}&-e_{n-2}\\
 \end{matrix}\right]}}}&~~~~A_{n-1}\oplus A_n\oplus E_{n-1} \ar[r]^-{\tiny{\left[\begin{matrix}
 \alpha_{n-1}&(-1)^{n}&0\\
 0&-g_{n}&-e_{n-1}\\
 \end{matrix}\right]}}\ar[d]^{\tiny{\left[\begin{matrix}
 1&0&0\\
 \alpha_{n-1}&(-1)^{n-2}&0\\
 g_{n-1}&0&(-1)^{n-1}\\
 \end{matrix}\right]}} & ~~~A_n\oplus E_n \ar[r]^-{\tiny{\left[\begin{matrix}
 \alpha_{n}&(-1)^{n}e_n\\
 \end{matrix}\right]}}\ar[d]^{\tiny{\left[\begin{matrix}
 1&0\\
 g_{n}&(-1)^{n}\\
 \end{matrix}\right]}} & A_{n+1} \ar@{-->}[r]^{\tiny{\left[\begin{matrix}
 -\alpha_0\\
 \varphi_0\\
 \end{matrix}\right]}_\ast\delta}\ar@{=}[d]^{~} &{~}\\
\cdots\ar[r]^-{\tiny{\left[\begin{matrix}
0&1&0\\
 0&0&0\\
0&0&e_{n-2}\\
 \end{matrix}\right]}} & A_{n-1}\oplus A_n\oplus E_{n-1} \ar[r]^{\tiny{\left[\begin{matrix}
0&1&0\\
 0&0&0\\
0&0&e_{n-1}\\
 \end{matrix}\right]}} & A_n\oplus E_n \ar[r]^{\tiny{\left[\begin{matrix}
0&e_{n}\\
 \end{matrix}\right]}} & A_{n+1} \ar@{-->}[r]^{{\varphi_0}_\ast\delta} & {~}}$$
shows that
$\textnormal{Coker}(\underline{k^\bullet})\cong I(\underline{\varphi^\bullet})$.

Similarly, we have
$$\textnormal{Ker}(\underline{c^{\bullet}}):=\xymatrix@C=2.3cm@R=2cm
{  B_0 \ar[r]^-{\tiny{\left[\begin{matrix}
e_0\\
 -\alpha_0\\
 \end{matrix}\right]}} & ~E_1\oplus B_1 \ar[r]^-{\tiny{\left[\begin{matrix}
-e_1&0\\
h_1&1\\
0&-\beta_1\\
 \end{matrix}\right]}}& E_2\oplus B_1\oplus B_2~~~\ar[r]^-{\tiny{\left[\begin{matrix}
-e_2&0&0\\
 h_2&\beta_1&1\\
0&0&-e_2\\
 \end{matrix}\right]}} &\cdots }$$
 $$\xymatrix@C=2.5cm@R=1.8cm{\cdots \ar[r]^-{\tiny{\left[\begin{matrix}
-e_{n-1}&0&0\\
 h_{n-1}&\beta_{n-2}&1\\
0&0&-\beta_{n-1}\\
 \end{matrix}\right]}} &~~E_{n}\oplus B_{n-1}\oplus B_{n}\ar[r]^{\tiny{\left[\begin{matrix}
-e_{n}&0&0\\
 h_{n}&\beta_{n-1}&1\\
 \end{matrix}\right]}} & A_{n+1}\oplus B_n \ar@{-->}[r]^{\tiny{\left[\begin{matrix}
 -\varphi_{n+1}&\beta_n\\
 \end{matrix}\right]}^{\ast}\delta'} & {~}}).$$
The following commutative diagram
 $$\xymatrix@C=2cm@R=1.8cm{
I(\varphi^\bullet)\ar[d]^-{\theta^\bullet} & B_0 \ar[r]^-{\tiny{\left[\begin{matrix}
 e_0\\
 0\\
 \end{matrix}\right]}} \ar@{=}[d] & E_1\oplus B_1 \ar[r]^-{\tiny{\left[\begin{matrix}
e_1&0\\
 0&1\\
 0&0\\
 \end{matrix}\right]}}\ar[d]^-{\tiny{\left[\begin{matrix}
 1&0\\
 -h_1&1\\
 \end{matrix}\right]}} & E_2\oplus B_1\oplus B_2\ar[r]^-{\tiny{\left[\begin{matrix}
 e_2&0&0\\
 0&0&1\\
 0&0&0\\
 \end{matrix}\right]}}\ar[d]^-{\tiny{\left[\begin{matrix}
 -1&0&0\\
 0&1&0\\
 h_2&-\beta_1&1\\
 \end{matrix}\right]}} & \cdots\\
\textnormal{Ker}(\underline{c^{\bullet}}) & B_0 \ar[r]^-{\tiny{\left[\begin{matrix}
 e_0\\
 -\beta_0\\
 \end{matrix}\right]}} &  E_1\oplus B_1  \ar[r]^-{\tiny{\left[\begin{matrix}
 -e_1&0\\
 h_1&1\\
 0&-\beta_1\\
 \end{matrix}\right]}} &E_2\oplus B_1\oplus B_2 \ar[r]^-{\tiny{\left[\begin{matrix}
  -e_2&0&0\\
 h_2&\beta_1&1\\
 0&0&-\beta_2\\
 \end{matrix}\right]}}& \cdots }$$
 $$\xymatrix@C=2.6cm@R=1.8cm{
 \cdots \ar[r]^-{\tiny{\left[\begin{matrix}
e_{n-1}&0&0\\
0&0&1\\
 0&0&0\\
 \end{matrix}\right]}}& ~E_n\oplus B_{n-1}\oplus B_n \ar[r]^-{\tiny{\left[\begin{matrix}
 e_{n}&0&0\\
0&0&1\\
 \end{matrix}\right]}}\ar[d]^{\tiny{\left[\begin{matrix}
 (-1)^{n-2}&0&0\\
 0&1&0\\
 (-1)^{n-1}h_{n-1}&-\beta_{n-2}&1\\
 \end{matrix}\right]}} & A_{n+1}\oplus B_n \ar@{-->}[r]^{{\varphi_0}_\ast\delta}\ar[d]^{\tiny{\left[\begin{matrix}
(-1)^{n}&0\\
0&1\\
 \end{matrix}\right]}} &{~}\\
\cdots\ar[r]^-{\tiny{\left[\begin{matrix}
-e_{n-1}&0&0\\
 h_{n-1}&\beta_{n-2}&1\\
0&0&-\beta_{n-1}\\
 \end{matrix}\right]}} &~~E_n\oplus B_{n-1}\oplus B_n\ar[r]^-{\tiny{\left[\begin{matrix}
-e_{n}&0&0\\
 h_{n}&\beta_{n-1}&1\\
 \end{matrix}\right]}} & A_{n+1}\oplus B_n \ar@{-->}[r]^-{\tiny{\left[\begin{matrix}
 -\varphi_{n+1}&\beta_n\\
 \end{matrix}\right]}^{\ast}\delta'} & {~}}$$
shows that
$I(\underline{\varphi^\bullet})\cong\textnormal{Ker}(\underline{c^\bullet})$. This concludes the proof \qed

\end{theorem}

\begin{remark}
According to Remark \ref{re43}, we know that when $\C$ is $n$-exact category, the conclusion Theorem \ref{th44} holds for both odd and even numbers of $n$.
\end{remark}

\begin{proposition}
Let $(\C,\mathbb{E},\mathfrak{s})$ be an $n$-exangulated category. Then a distinguished $n$-exangle
$$P_A^\bullet=(\xymatrix{
 TA\ar[r]^-{\alpha_0} & P_1^A \ar[r]^-{\alpha_1}& P_2^A \ar[r]^-{\alpha_2} & \cdots \ar[r]^-{\alpha_{n-2}}&P_{n-1}^A \ar[r]^-{\alpha_{n-1}}& P_n^A \ar[r]^-{\alpha_{n}} & A \ar@{-->}[r]^-{{\rho}^A}&})$$
with $P_1^A,P_2^A,\cdots,P_n^A\in\P$ is a projective object in $S(\C)/\R_2$. Moreover, if $\C$ has enough projectives, then each projective object in $S(\C)/\R_2$ is of the form $P_A^\bullet$.

\proof Suppose that $\underline{\varphi^\bullet}:B^\bullet\rightarrow C^\bullet$ is an epimorphism
$$\xymatrix{
B^{\bullet}\ar[d]^-{\varphi^\bullet} & B_0 \ar[r]^{\beta_0}\ar[d]^{\varphi_0} & B_1 \ar[r]^{\beta_1}\ar[d]^{\varphi_1} & B_2 \ar[r]^{\beta_2}\ar[d]^{\varphi_2} & \cdots \ar[r]^{\beta_{n-1}}& B_{n} \ar[r]^{\beta_{n}}\ar[d]^{\varphi_{n}} & B_{n+1} \ar@{-->}[r]^{\delta}\ar[d]^{\varphi_{n+1}} &{~}\\
C^{\bullet} & C_0 \ar[r]^{\gamma_0} & C_1 \ar[r]^{\gamma_1} & C_2 \ar[r]^{\gamma_2} & \cdots \ar[r]^{\gamma_{n-1}}& C_{n} \ar[r]^{\gamma_{n}} & C_{n+1} \ar@{-->}[r]^{\delta'} & {~}}$$
and $\underline{\psi^\bullet}:P_A^\bullet\rightarrow C^\bullet$ is a morphism in $S(\C)/\R_2$.
$$\xymatrix{
P_A^{\bullet}\ar[d]^-{\psi^\bullet} & TA \ar[r]^{\beta_0}\ar[d]^{\psi_0} & P_1^A \ar[r]^{\beta_1}\ar[d]^{\psi_1} & P_2^A \ar[r]^{\beta_2}\ar[d]^{\psi_2} & \cdots \ar[r]^{\beta_{n-1}}& P_{n}^A \ar[r]^{\beta_{n}}\ar[d]^{\varphi_{n}} & A \ar@{-->}[r]^{\rho^A}\ar[d]^{\psi_{n+1}} &{~}\\
C^{\bullet} & C_0 \ar[r]^{\gamma_0} & C_1 \ar[r]^{\gamma_1} & C_2 \ar[r]^{\gamma_2} & \cdots \ar[r]^{\gamma_{n-1}}& C_{n} \ar[r]^{\gamma_{n}} & C_{n+1} \ar@{-->}[r]^{\delta'} & {~}}$$
According to Lemma \ref{lemma42}(3), we have the morphism ${\left[\begin{matrix}
 \varphi_{n+1}&\gamma_n\\
 \end{matrix}\right]}:B_{n+1}\oplus C_n\rightarrow C_{n+1}$ is a retraction. Hence, there exists a morphism ${\left[\begin{matrix}
 \varphi_{n+1}'\\
 \gamma_n'
 \end{matrix}\right]}:C_{n+1}\rightarrow B_{n+1}\oplus C_n$ such that $\varphi_{n+1}\varphi_{n+1}'+\gamma_n\gamma_n'=1.$  We set $\phi_{n+1}=\varphi_{n+1}'\psi_{n+1}:A\rightarrow B_{n+1}$. By Lemma \ref{lemma212}, we have the following exact sequence:
 $$\C(TA, B_0) \xrightarrow{(\rho^A)^\sharp} \mathbb{E}(A, B_0)\xrightarrow{(\alpha_n)^\ast} \mathbb{E}(P_n^A, B_0) $$
  Since, the object $P_n^A$ is a projective, by Lemma \ref{lemma213}, $\mathbb{E}(P_n^A, B_0)=0$, there exists a morphism $\phi_0:TA\rightarrow B_0$ such that ${\phi_0}_\ast\rho^A=(\phi_{n+1})^\ast\delta$. This gives a morphism of extensions $(\phi_0,\phi_{n+1}):\rho^A\rightarrow \sigma$. That is
  $$\xymatrix{
P_A^{\bullet} \ar@{-->}[d]^-{\phi^\bullet} &TA \ar[r]^{\beta_0} \ar@{-->}[d]^{\phi_0} & P_1^A \ar[r]^{\beta_1} \ar@{-->}[d]^{\phi_1} & P_2^A \ar[r]^{\beta_2} \ar@{-->}[d]^{\phi_2} & \cdots \ar[r]^{\beta_{n-1}}& P_{n}^A \ar[r]^{\beta_{n}} \ar@{-->}[d]^{\phi_{n}} & A \ar@{-->}[r]^{\rho^A}\ar[d]^{\phi_{n+1}} &{~}\\
B^{\bullet} & B_0 \ar[r]^{\beta_0} & B_1 \ar[r]^{\beta_1} & B_2 \ar[r]^{\beta_2} & \cdots \ar[r]^{\beta_{n-1}}& B_{n} \ar[r]^{\beta_{n}} & B_{n+1} \ar@{-->}[r]^{\delta} & {.}} $$
 Note that $$\psi_{n+1}-\varphi_{n+1}\phi_{n+1}=\psi_{n+1}-\varphi_{n+1}\varphi_{n+1}'\psi_{n+1}
 =(1-\varphi_{n+1}\varphi_{n+1}')=\gamma_n(\gamma_n'\psi_{n+1}),$$
 by Lemma \ref{lemma42}(1), we obtain $\underline{\psi^\bullet}=\underline{\varphi^\bullet\phi^\bullet}$. Therefore, $P_A^\bullet$ is a  projective object in $S(\C)/\R_2$.

 Assume that $\C$ has enough projectives, For each distinguished $n$-exangle
 $$A^\bullet=(\xymatrix{
 A_0\ar[r]^{\alpha_0} & A_1 \ar[r]^{\alpha_1}& A_2 \ar[r]^{\alpha_2} & \cdots \ar[r]^{\alpha_{n-2}}&A_{n-1} \ar[r]^{\alpha_{n-1}}& A_n \ar[r]^{\alpha_{n}} & A_{n+1} \ar@{-->}[r]^{\delta}&})$$
there exists a distinguished $n$-exangle
 $$P_{A_{n+1}}^\bullet=(\xymatrix{
 TA_{n+1}\ar[r]^-{\alpha_0'} & P_1^{A_{n+1}}\ar[r]^-{\alpha_1'}& \cdots \ar[r]^-{\alpha_{n-2}'}&{P_{n-1}}^{A_{n+1}} \ar[r]^-{\alpha_{n-1}'}& {P_{n}}^{A_{n+1}} \ar[r]^-{\alpha_{n}'} & A_{n+1} \ar@{-->}[r]^{\rho^{A_{n+1}}}&})$$
 with $P_1^{A_{n+1}},P_2^{A_{n+1}},\cdots,{P_{n}}^{A_{n+1}}\in\P$. By the dual of Lemma \ref{lemma215}, we have a morphism
 $\varphi^\bullet:P_{A_{n+1}}^\bullet\rightarrow A^\bullet$ as follows:
 $$\xymatrix{
P_{A_{n+1}}^{\bullet} \ar@{-->}[d]^-{\varphi^\bullet} & TA_{n+1} \ar[r]^{\alpha_0'} \ar@{-->}[d]^{\varphi_0} &P_1^{A_{n+1}} \ar[r]^{\alpha_1'} \ar@{-->}[d]^{\varphi_1} & P_2^{A_{n+1}} \ar[r]^{\alpha_2'} \ar@{-->}[d]^{\varphi_2} & \cdots \ar[r]^{\alpha_{n-1}'}& {P_{n}}^{A_{n+1}} \ar[r]^{\alpha_{n}'} \ar@{-->}[d]^{\varphi_{n}} & A_{n+1} \ar@{-->}[r]^{\rho^{A_{n+1}}}\ar@{=}[d]^{~} &{~}\\
A^{\bullet} & A_0 \ar[r]^{\alpha_0} & A_1 \ar[r]^{\alpha_1} & A_2 \ar[r]^{\alpha_2} & \cdots \ar[r]^{\alpha_{n-1}}& A_{n} \ar[r]^{\alpha_{n}} & A_{n+1} \ar@{-->}[r]^{\delta} & {.}} $$
Since ${\left[\begin{matrix}
 1&0\\
 \end{matrix}\right]}{\left[\begin{matrix}
 1\\
\alpha_n
 \end{matrix}\right]}=1$, hence $\underline{\varphi^\bullet}:P_{A_{n+1}}^\bullet\rightarrow A^\bullet$ is an epimorphism and $P_{A_{n+1}}^\bullet$ is projective. Furthermore, suppose that $A^\bullet$ is projective object in $S(\C)/\R_2$, then the epimorphism $\underline{\varphi^\bullet}:P_{A_{n+1}}^\bullet\rightarrow A^\bullet$ is split. Therefore each projective object of $S(\C)/\R_2$ is of the form $P_A^\bullet$ for some object $A$ in
 $\C$. \qed
\end{proposition}

\begin{definition}\textnormal{\cite[Definition 3.1]{HHZ}}\label{df47}
Let $\C$ be an $n$-exangulated category. A distinguished $n$-exangle
 $$\xymatrix{
 A_0\ar[r]^{\alpha_0} & A_1 \ar[r]^{\alpha_1}& A_2 \ar[r]^{\alpha_2}& \cdots \ar[r]^{\alpha_{n-2}}&A_{n-1} \ar[r]^{\alpha_{n-1}}& A_n \ar[r]^{\alpha_{n}} & A_{n+1} \ar@{-->}[r]^{\delta}&}$$
in $\C$ is called an \emph{Auslander-Reiten $n$-exangle} if  the following holds:
\begin{itemize}
\item[$\rm{(1)}$] $\delta\in\mathbb{E}(A_{n+1},A_0)$ is non split.
\item[$\rm{(2)}$] If $\beta:A_0\rightarrow B$ is not a section, then $\beta$ factors through $\alpha_0$.
\item[$\rm{(3)}$] If $\gamma:C\rightarrow A_{n+1}$ is not a retraction, then $\gamma$ factors through $\alpha_{n}$.
\item[$\rm{(4)}$] when $n\geq 2, \alpha_1,\alpha_2,\cdots,\alpha_{n-1}$ are in ${\textnormal{rad}_{\C}}$.
\end{itemize}
\end{definition}

\begin{remark}
\begin{itemize}
\item[$\rm{(1)}$] If $\C$ is an $n$-abelian category, then Definition \ref{df47} coincides with the definition of $n$-Auslander-Reiten sequence of $n$-abelian category (cf.\cite{I,XL}), which is first introduced by Iyama in \cite[Definition 3.1]{I} .
\item[$\rm{(2)}$] If $\C$ is an $(n+2)$-angulated category, then Definition \ref{df47} coincides with the definition of Auslander-Reiten $(n+2)$-angle of $(n+2)$-angulated category (cf.\cite{F,Z}).
\end{itemize}
\end{remark}

\begin{proposition}
Let $\C$ be a Krull-Smidt $n$-exangulated category. when taking $n$ as an even number, suppose that $A^\bullet: \xymatrix{
 A_0\ar[r]^{\alpha_0} & A_1 \ar[r]^{\alpha_1}& A_2 \ar[r]^{\alpha_2}& \cdots \ar[r]^{\alpha_{n-2}}&A_{n-1} \ar[r]^{\alpha_{n-1}}& A_n \ar[r]^{\alpha_{n}} & A_{n+1} \ar@{-->}[r]^{\delta}&}$
 is a non-split distinguished $n$-exangle such that $A_0$ and $A_{n+1}$ are indecomposable and $\alpha_1,\alpha_2,\cdots,\alpha_{n-1}\in{\textnormal{rad}_{\C}}$. Then $A^\bullet$ is simple object in $S(\C)/\R_2$ if and only if $A^\bullet$ is an Auslander-Reiten distinguished $n$-exangle in $\C$.
 \proof Suppose that $A^\bullet$ is simple object in $S(\C)/\R_2$, we need to proof  $A^\bullet$ is an Auslander-Reiten distinguished $n$-exangle in $\C$. Assume that $\varphi_0:A_0\rightarrow B_0$ is not a section. By the dual of Lemma \ref{lemma215}, we have a morphism $\varphi^\bullet:A^{\bullet}\rightarrow B^\bullet$ as follows:
 $$\xymatrix{
A^{\bullet} \ar[d]^-{\varphi^\bullet}&A_0 \ar[r]^{\alpha_0}\ar[d]^{\varphi_0} & A_1 \ar[r]^{\alpha_1}\ar[d]^{\varphi_1} & A_2 \ar[r]^{\alpha_2}\ar[d]^{\varphi_2} & \cdots \ar[r]^{\alpha_{n-1}}& A_{n} \ar[r]^{\alpha_{n}}\ar[d]^{\varphi_{n}} & A_{n+1} \ar@{-->}[r]^{\delta}\ar@{=}[d]^{~} &{~}\\
B^\bullet & B_0 \ar[r]^{\beta_0} & B_1 \ar[r]^{\beta_1} & B_2 \ar[r]^{\beta_2} & \cdots \ar[r]^{\beta_{n-1}}& B_{n} \ar[r]^{\beta_{n}} & B_{n+1} \ar@{-->}[r]^{{\varphi_0}_\ast\delta} & {~}}$$

Since $\alpha_0$ and $\varphi_0$ are not sections and $A_0$  is indecomposable, hence, the morphism ${\left[\begin{matrix}
\alpha_0\\
\varphi_0\\
 \end{matrix}\right]}:A_0\rightarrow A_1\oplus B_0$ is not a section. Therefore, by the dual of Lemma \ref{lemma42}(3), $\underline{\varphi^\bullet}$ is not a monomorphism. Since $A^\bullet$ is simple object, thus, $\underline{\varphi^\bullet}=0$. It follow from Lemma \ref{lemma42}(1) that $\varphi_0$ factor through $\alpha_0$. Dually, we can prove that if $\varphi_{n+1}:C\rightarrow A_{n+1}$ is not a retraction, then $\varphi_{n+1}$ factors through $\alpha_n$. Therefore, $A^\bullet$ is an Auslander-Reiten distinguished $n$-exangle in $\C$.

 Conversely, suppose that $\varphi^\bullet:A^\bullet\rightarrow B^\bullet$ is a morphism of distinguished $n$-exangle. On the one hand, if $\varphi_0$ is a section, obviously, $\left[\begin{matrix}
\alpha_0\\
\varphi_0\\
 \end{matrix}\right]:A_0\rightarrow A_1\oplus B_0$ is also. It follows from Lemma \ref{lemma42}(3) that $\underline{\varphi^\bullet}$ is a monomorphism. On the other hand, if $\varphi_0$ is not a section, since $A^\bullet$ is an Auslander-Reiten distinguished $n$-exangle, therefore, $\varphi_0$ factors through $\alpha_0$, thus $\underline{\varphi^\bullet}=0$. Hence, each morphism $\underline{\varphi^\bullet}:A^\bullet\rightarrow B^\bullet$ either a monomorphism or a zero morphism. That is  $A^\bullet$ is simple object in $S(\C)/\R_2$. This concludes the proof. \qed
\end{proposition}

From now on, we assume that $\C$ is an $n$-exangulated category with enough projectives $\P$ and injectives $\I$.

\begin{definition}
Let $\C$ be an $n$-exangulated category. Given a distinguished $n$-exangle $$\delta=(\xymatrix{
 A_0\ar[r]^{\alpha_0} & A_1 \ar[r]^{\alpha_1}& A_2 \ar[r]^{\alpha_2}& \cdots\ar[r]^{\alpha_{n-2}}& A_{n-1}  \ar[r]^{\alpha_{n-1}}& A_n \ar[r]^{\alpha_{n}} & A_{n+1} \ar@{-->}[r]^{\rho}&}),$$
 $\delta^\ast$  is called a \emph{contravariant defect} and $\delta_\ast$  is called a \emph{convariant defect} if the following sequences of functors are exact
 $$\C(-, A_0) \xrightarrow{\C(-,\alpha_0)} \C(-, A_1)\xrightarrow{\C(-,\alpha_1)}\cdots\xrightarrow{\C(-,\alpha_{n-1})}
  \C(-,A_n)\xrightarrow{\C(-,\alpha_n)} \C(-,A_{n+1})\rightarrow\delta^\ast\rightarrow0,$$
 $$\C(A_{n+1},-) \xrightarrow{\C(\alpha_n,-)} \C(A_{n},-)\xrightarrow{\C(\alpha_{n-1},-)}\cdots\xrightarrow{\C(\alpha_{1},-)}
  \C(A_1,-)\xrightarrow{\C(\alpha_0,-)} \C(A_0,-)\rightarrow\delta_\ast\rightarrow0.$$

\end{definition}

\begin{example}\label{ex411}
\begin{itemize}
\item[$\rm{(1)}$] Let $\delta=P_A^\bullet=(\xymatrix{
 TA\ar[r]^-{\alpha_0} & P_1^A \ar[r]^-{\alpha_1}& \cdots \ar[r]^-{\alpha_{n-2}}&P_{n-1}^A \ar[r]^-{\alpha_{n-1}}& P_n^A \ar[r]^-{\alpha_{n}} & A \ar@{-->}[r]^-{{\rho}}&})$ with $P_1^{A_{n+1}},P_2^{A_{n+1}},\cdots,{P_{n}}^{A_{n+1}}\in\P$. Then $\delta^\ast=(\C/[\P])(-,A)$ and $\delta_\ast=\mathbb{E}(A,-)$.

 \item[$\rm{(2)}$] Let $\delta=I_A^\bullet=(\xymatrix{
 A\ar[r]^-{\alpha_0} & I_1^A \ar[r]^-{\alpha_1}&I_2^A \ar[r]^-{\alpha_1}& \cdots \ar[r]^-{\alpha_{n-2}}&I_{n-1}^A \ar[r]^-{\alpha_{n-1}}& I_n^A \ar[r]^-{\alpha_{n}} & SA \ar@{-->}[r]^-{{\rho}}&})$ with \\ $I_1^{A_{n+1}},I_2^{A_{n+1}},\cdots,{I_{n}}^{A_{n+1}}\in\I$. Then $\delta^\ast=\mathbb{E}(-,A)$ and $\delta_\ast=(\C/[\I])(A,-).$
\end{itemize}

\end{example}

\begin{theorem}\label{th412}
Let $(\C,\mathbb{E},\mathfrak{s})$ be an $n$-exangulated category with enough projectives $\P$ and enough injectives $\I$.
\begin{itemize}
\item[$\rm{(1)}$] We have the following equivalences
$$S(\C)/\R_2\cong\textnormal{mod-}(\C/[\P])\cong (\textnormal{mod-}(\C/[\I])^{op})^{op}.$$
Furthermore, the equivalence $F:S(\C)/\R_2\cong\textnormal{mod-}(\C/[\P])$ is given by $\delta\mapsto\delta^\ast$ and the equivalence
$G: S(\C)/\R_2\cong(\textnormal{mod-}(\C/[\I])^{op})^{op}$ is given by $\delta\mapsto\delta_\ast.$
\item[$\rm{(2)}$] The abelian category $\textnormal{mod-}(\C/[\P])$ has enough projectives and enough injectives. Furthermore, each projective object is of the form $(\C/[\P])(-,A)$,  and when taking $n$ as an even number, each injective object is of the form $\mathbb{E}(-,A)$.
\item[$\rm{(3)}$] The abelian category $\textnormal{mod-}(\C/[\I])^{op}$ has enough projectives and enough injectives. Furthermore, each projective object is of the form $(\C/[\I])(A,-)$, and when taking $n$ as an even number, each injective object is of the form $\mathbb{E}(A,-)$.
\end{itemize}
\proof (1) The first assertion follows directly from Theorem \ref{th41}, Assume that $$\delta=(\xymatrix{
 A_0\ar[r]^{\alpha_0} & A_1 \ar[r]^{\alpha_1}& A_2 \ar[r]^{\alpha_2}& \cdots\ar[r]^{\alpha_{n-2}}& A_{n-1}  \ar[r]^{\alpha_{n-1}}& A_n \ar[r]^{\alpha_{n}} & A_{n+1} \ar@{-->}[r]^{\rho}&})$$
 is a distinguished $n$-exangle. According to the definition of $\delta^\ast$, we have $\delta^\ast=\textnormal{Coker}(\C(-,\alpha_n))$. Since $\delta^\ast(\P)=0$, hence, we define $F(\delta)=\textnormal{Coker}(\C/[\P](-,\alpha_n))$, we view $\delta^\ast$ as a finitely presented $(\C/[\P])$-module. Therefore, $F(\delta)=\delta^\ast$. Similarly, we have $G(\delta)=\delta_\ast$.

 (2) and (3) follows from Theorem \ref{th412}(1), Theorem \ref{th44} and Example
 \ref{ex411}. \qed

\end{theorem}

\begin{corollary}\textnormal{\cite[Corollary 4.9]{Li}}
Let $(\C,\mathbb{E},\mathfrak{s})$ be an $n$-exangulated category with enough projectives $\P$ and enough injectives $\I$. Then there is a duality
$$\Phi:\textnormal{mod-}(\C/[\P])\rightarrow\textnormal{mod-}(\C/[\I])^{op},~~~
\delta^\ast\mapsto\delta_\ast.$$
Furthermore, by taking $n$ as an even number, we obtain the following two dulities
$$\Phi:\textnormal{proj-}(\C/[\P])\rightarrow\textnormal{inj-}(\C/[\I])^{op},~~~
(\C/[\P])(-,A)\mapsto\mathbb{E}(A,-).$$
$$\Phi:\textnormal{inj-}(\C/[\P])\rightarrow\textnormal{proj-}(\C/[\I])^{op},~~~
\mathbb{E}(-,A)\mapsto(\C/[\I])(A,-).$$
\proof It follows directly from Theorem \ref{th412}. \qed
\end{corollary}

\textbf{Yutong Zhou}\\
School of Mathematics and statistics, Northeast Normal University, 130024 Changchun, Jilin, P. R. China.\\
E-mail: \textsf{yutongzhou2021@163.com}
\\[0.3cm]

\end{document}